\documentclass{article}
\usepackage{graphicx}
\usepackage{amsmath}
\usepackage{amssymb}

\newcommand{\RR}{\mathbb{R}}

\newcommand{\CC}{\mathbb{C}}

\newcommand{\nocode}[1]{}

  \newtheorem{prop}{Proposition}

  \newtheorem{rmk}{Remark}

  \newtheorem{alg}{Algorithm}
  \newtheorem{method}{Method}
  \newtheorem{conj}{Conjecture}
  \newtheorem{quest}{Question}

\begin{document}
\thispagestyle{empty}

\begin{center}
\begin{Large}
{\bf Solutions of Polynomial Systems Derived from the Steady Cavity Flow Problem}
\end{Large}
 \vspace{5mm}\\ 
$\mbox{Martin Mevissen}^{\star}$,
$\mbox{Kosuke Yokoyama}^{\dagger}$,
$\mbox{Nobuki Takayama}^{\sharp}$  \vspace{5mm} \\
November 14, 2008
\end{center} 
\vspace{1cm}

\noindent
{\bf Abstract.} \\ 
We propose a general algorithm to enumerate all solutions of a zero-dimensional polynomial system with respect to a given cost function. 
The algorithm is developed and is used to study a polynomial system obtained by discretizing the steady cavity flow problem in two dimensions. 
The key technique on which our algorithm is based  is to solve polynomial optimization problems via sparse semidefinite programming relaxations (SDPR) \cite{wkkm}, which has been adopted successfully to solve reaction-diffusion boundary value problems in \cite{mknt}. 
The cost function to be minimized is derived from discretizing the fluid's kinetic energy. The enumeration algorithm's solutions are shown to converge to the minimal kinetic energy solutions for SDPR of increasing order. We demonstrate the algorithm with SDPR of first and second order on polynomial systems for different scenarios of the cavity flow problem and succeed in deriving the $k$ smallest kinetic energy solutions.  The question whether these solutions converge to solutions of the steady cavity flow problem is discussed, and we pose a conjecture for the minimal energy solution for increasing Reynolds number.
\vspace{1cm}
  
\noindent
{\bf Key words. } \vspace{0.1cm} \\
Steady cavity flow problem, finite difference discretization,
polynomial optimization, 
semidefinite programming relaxation, sparsity
\vspace{1cm}

\noindent
\parbox[t]{0.5cm}{$\star$}
\parbox[t]{14.9cm}{Department of Mathematical and Computing Sciences, Tokyo Institute of Technology, 2-12-1 Ookayama, Meguro-ku, Tokyo 152-8552 Japan. 
{\it martime6@is.titech.ac.jp}. Research supported by the Doctoral Scholarship of the German Academic Exchange Service (DAAD).}
\medskip

\noindent
\parbox[t]{0.5cm}{$\dagger$}
\parbox[t]{14.9cm}{Department of Mathematics, Graduate School of Science, Kobe University, 
1-1 Rokkodai, Nada-ku, Kobe 657-8501 Japan. Current affiliation is Fujitsu Cooporation}
\medskip

\noindent
\parbox[t]{0.5cm}{$\sharp$}
\parbox[t]{14.9cm}{Department of Mathematics, Graduate School of Science, Kobe University, 
1-1 Rokkodai, Nada-ku, Kobe 657-8501 Japan. 
{\it  takayama@math.kobe-u.ac.jp}
}
\newpage
\pagenumbering{arabic}

\section{Introduction}
The steady cavity flow problem is a simple model of a flow with closed streamlines and is used for examining and validating numerical solution techniques in fluid dynamics. Although it has been discussed in several literature of numerical analysis of fluid mechanics
(see, e.g., \cite{kawaguti}, \cite{burggraf}, \cite{gustafson}, \cite{cheng}, \cite{takami}),
it is still an interesting problem to a number of researchers for a range of Reynolds numbers.
We are interested in a polynomial system derived from discretizing the steady cavity flow problem.
This polynomial system, called the {\it discrete cavity flow problem}, is obtained by discretizing the cavity region, approximating the partial differential equation of the two-dimensional cavity flow problem by finite difference method, and depends on two parameters, the Reynolds number $R$ and the boundary velocity $v$. 

Our main contribution is an algorithm to enumerate the discrete cavity flow problem's solutions with respect to an objective function, that is derived from discretizing the cavity flow's kinetic energy function. The key element of the enumeration algorithm is the sparse semidefinite programming relaxation method (SDPR) \cite{wkkm} for solving polynomial optimization problems, whose solution is taken as the starting point for Newton's method or sequential quadratic programming. Recently, the SDPR has been successfully adopted to derive numerical solutions to a class of reaction diffusion equations \cite{mknt}.  In this paper, the polynomial optimization problem is the minimization of the discretized kinetic energy subjected to the discrete cavity flow problem. 
We prove that the first $k$ solutions provided by the enumeration algorithm converge to the $k$ smallest energy solutions of the discrete cavity flow problem, in case that we apply SDPR of increasing relaxation order. Furthermore, we demonstrate this algorithm for different parameter settings of $R$ and $v$, and show in some examples that it is sufficient to apply SDPR with first or second order, to enumerate accurate approximations to the smallest energy solutions.
At second, we discuss the minimal energy solution's behavior of the discrete steady cavity flow problem in case that a finer grid is chosen to discretize the cavity flow problem. For small Reynolds numbers $R$ standard grid-refining techniques can be  applied to extend solutions of the polynomial system to finer grids. 
However the polynomial systems for large $R$ and $v$ behave differently and convergence is far more difficult to obtain. We examine the polynomial systems for a fixed discretization and increasing Reynolds number $R$.
Based on our observations, we conjecture the minimal kinetic energy is converging to zero if $R$ tends to infinity. Also, we test the performance of SDPR for an alternative finite difference discretization by Arakawa \cite{arakawa} of the steady cavity flow problem.\\
The exact formulation of the steady cavity flow problem in 2 dimensions is introduced in section 2. We discuss its boundary conditions and derive the discrete steady cavity flow problem. In section 3 we study the discrete cavity flow problem by Gr\"obner basis method for coarse grid discretizations, in order to be able to verify the results derived by our enumeration algorithm later on. 
In section 4 we show how to solve a polynomial optimization problem derived from the discrete cavity flow problem via the SDPR method. We present our main contribution, the algorithm to enumerate the discrete cavity flow problems solutions w.r.t. their kinetic energy, and demonstrate its power for some settings of the 2 parameters $R$ and $v$.
In section 5 we discuss the question of how to extend the discrete cavity flow problem's minimal energy solution to finer grids and its convergence to analytic solutions of the steady cavity flow problem. Finally, we examine the discrete steady cavity flow problem for increasing Reynolds number $R$.

\section{A Polynomial System for the Steady Cavity Flow Problem}

Let us review the well-known stream function method to solve the Navier-Stokes
equation 
(see, e.g., \cite{kawaguti}, \cite{takami}).
The stream function method is a standard method to solve the 2 dimensional
steady cavity flow problem numerically. 

Let $(u(x,y,t),v(x,y,t))$ be the velocity of the two dimensional cavity
flow of an incompressible fluid.
It follows from the continuity equation of the incompressible fluid
(preservation of the mass)
  $\frac{\partial u}{\partial x} + \frac{\partial v}{\partial y} = 0$
that
there exists a function $\psi(x,y,t)$ such that
\begin{equation}
  \frac{\partial \psi}{\partial x} = -v, \ 
  \frac{\partial \psi}{\partial y} = u.
\end{equation}
Put ${\bf v} = (u,v,0)$.
${\vec \omega} = {\rm rot}\, {\bf v}$ is called the vorticity.
Since the last coordinate of ${\bf v}$ is $0$,
${\vec \omega}$ can be written as $(0,0,\omega(x,y,t))$.
The continuity equation and the Navier-Stokes equation (preservation 
of the momentum) can be written as follows in terms of $\psi$ and $\omega$.
\begin{eqnarray}
 & & \Delta \psi = - \omega   \label{eq:st1} \\
 & & \frac{\partial \omega}{\partial t} =
    \frac{\partial \psi}{\partial y}
    \frac{\partial \omega}{\partial x} - 
    \frac{\partial \psi}{\partial x}
    \frac{\partial \omega}{\partial y} + 
    \frac{1}{R} \Delta \omega  \label{eq:st2}
\end{eqnarray}
Here, $\Delta$ is the Laplace operator and $R$ is the Reynolds number.
Let us consider the cavity region $ABCD$ with the coordinate
$A=(0,0)$, $B=(0,-1)$, $C=(1,-1)$, $D=(1,0)$.

\begin{center}
\setlength{\unitlength}{0.2mm}
\begin{picture}(130,130)(-5,-105)
\put(40,0){\vector(1,0){20}}
\put(0,0){\line(0,-1){100}}
\put(0,-100){\line(1,0){100}}
\put(100,-100){\line(0,1){100}}
\put(0,0){$A$}
\put(0,-100){$B$}
\put(100,-100){$C$}
\put(100,0){$D$}
\end{picture}
\end{center}

The steady cavity flow problem is (\ref{eq:st1}) and
(\ref{eq:st2}) with the steady condition 
$\frac{\partial \omega}{\partial t}=0$
and the boundary condition
\begin{eqnarray}
 && u(0,y) = u(x,0) = u(1,y) = 0 \quad \mbox{ on } AB, BC, CD \label{eq:b1} \\
 && v(0,y) = v(x,0) = v(1,y) = 0 \quad \mbox{ on } AB, BC, CD \label{eq:b2}\\
 && u(x,1) = s, v(x,1) = 0 \quad \mbox{ on } AD \label{eq:b3}
\end{eqnarray}
Here $s$ is the velocity of the stream out of the cavity $ABCD$. 

We devide the square $ABCD$ into the $N \times N$ mesh. 
Put $h = 1/N$.
Let us translate these boundary conditions into boundary conditions for
$\psi$ and $\omega$.
It follows from (\ref{eq:b1}), (\ref{eq:b2}), (\ref{eq:b3}) that
the function $\psi$ is constant on the boundaries $AB$, $BC$, $CD$, $DA$.
Since $\psi$ is continuous, we suppose that $\psi=0$ on the boundaries. 

The boundary condition for $\omega$ is a little complicated; 
see, e.g.,  \cite[p.162]{takami}.
We cite this discussion in text books for reader's convinience.
Let us consider the case of the boundary $AD$.
We take a mesh point $M$ on $AD$.
Let $P$ be the mesh point inside the cavity ajacent to $M$
and $P'$ the mirror image of $P$ with respect to $AD$.
We supposed that the size of the mesh is $h$.

\setlength{\unitlength}{1mm}
\begin{picture}(130,40)(-5,-20)
\put(0,0){\line(1,0){100}}
\put(37,-10){\line(1,0){26}}
\put(37,10){\line(1,0){26}}
\put(40,-13){\line(0,1){26}}
\put(50,-13){\line(0,1){26}}
\put(60,-13){\line(0,1){26}}
\put(0,0){$A$}
\put(100,0){$D$}
\put(50,-10){$P$}
\put(50,10){$P'$}
\put(46,-5){$h$}
\put(50,0){$M$}
\end{picture}

We denote the value of $\psi$ at the point $P$ by
$\psi(P)$ or $\psi_P$.
We have
$-\omega(M)=\Delta \psi (M) = \psi_{yy} \simeq 
\frac{\psi_{P} - 2 \psi_M + \psi_{P'}}{h^2}$.
We need to determine the value of $\psi_{P'}$ to
get an approximate value of $\omega$ at $M$.
By using the central difference approximation,
$s=v=\frac{\partial \psi}{\partial y}(M) \simeq \frac{\psi_{P'}-\psi_P}{2h}$.
Then, $\psi_{P'} \simeq 2hs + \psi_P$.
Therefore, we have
\begin{equation} \label{eq:boundary1}
  \omega_M \simeq -\frac{2 \psi_{P}+ 2hs}{h^2}
\end{equation}
Analogously, we have
\begin{equation} \label{eq:boundary2}
  \omega_M \simeq -\frac{2 \psi_{P}}{h^2}
\end{equation} 
when $M$ is a grid point on $AB$ or $BC$ or $CD$ and
$P$ is the adjacent internal grid point of $M$.

It follows from the discussion above, we obtain the following central difference
scheme for the steady cavity flow problem.
{\footnotesize
\begin{eqnarray}
\nonumber \\
  0  &=&  -4\omega_{i,j}+\omega_{i+1,j}+\omega_{i-1,j}+\omega_{i,j+1}+\omega_{i,j-1} \nonumber\\
&&\ \ \ \ \ +\dfrac{R}{4}\{ (\psi_{i+1,j}-\psi_{i-1,j})(\omega_{i,j+1}-\omega_{i,j-1}) \label{eq:omegaconst}\\
&&\ \ \ \ \ \ \ \ \ \ \ \ -(\psi_{i,j+1}-\psi_{i,j-1})(\omega_{i+1,j}-\omega_{i-1,j}) \nonumber \}
\\
0 &=& -4\psi_{i,j} +\psi_{i+1,j}+\psi_{i-1,j}+\psi_{i,j+1}+\psi_{i,j-1}+h^2\omega_{i,j}
\label{eq:psiconst}
\end{eqnarray}
}
$\psi = 0$ on the boundaries and
\begin{equation}  \label{eq:boundary}
\left\{
\begin{array}{c c c}
\omega = -\dfrac{2\psi_{P}+2sh}{h^2} &  (\mbox{on } AD)\\
\omega = -2\dfrac{\psi_P}{h^2} & (\mbox{on } AB,BC,CD) 
\end{array}
\right.
\end{equation}
We will call the polynomial system 
(\ref{eq:omegaconst}),
(\ref{eq:psiconst}),
(\ref{eq:boundary})
the {\it discrete steady cavity flow problem} denoted as DSCF($R,s,N$). Let the number $2(N-2)^2$ of variables in DSCF($R,s,N$) corresponding to function evaluations at interior grid points, be called the {\bf dimension $n$ of the discrete cavity flow problem}. Moreover, a solution $(\psi,\omega)(N)$ of the discrete cavity flow problem of discretization $N$, that does not converge to a physical solution of the original continuous cavity flow problem for $N\rightarrow\infty$ is called a {\bf fake solution}. 

Several methods have been used to solve the cavity flow problem and
the steady cavity flow problem numerically
(see, e.g., \cite{arakawa}, \cite{burggraf}, \cite{cheng}, \cite{gustafson}
\cite{kawaguti}, \cite{takami}, \cite{strikwerda}).
In this paper, we propose a new method to solve the discrete steady 
cavity flow problem.
This method provides solutions sorted by their (discretized) kinetic energy.

\begin{rmk}
\label{rmkAra}
In (\ref{eq:omegaconst}), we discretize the Jacobian
$ \frac{\partial \psi}{\partial y} \frac{\partial \omega}{\partial x} -
  \frac{\partial \psi}{\partial x} \frac{\partial \omega}{\partial y}$
by the central difference scheme. 
It is shown by Arakawa \cite{arakawa} that the central difference scheme is the simplest,
but the discretized system does not keep important physical invariants.
We study the system DSCF($R,s,N$) as the simplest starting test case.
\end{rmk}

\begin{rmk} 
We conjecture that the discrete cavity flow problem DSCF($R,s,N$) has finite complex solutions.
In other words, it defines a zero-dimensional ideal.  
We have checked the conjecture up to $N=5$ by Gr\"obner basis computation. 
\end{rmk}
In the sequel, the boundary velocity $s$ will be denoted by $v$ as long as no confusion
arises with $v=-\partial \psi/\partial x$.

\section{Gr\"obner Basis Method}   \label{section:gbmethod}

The Gr\"obner basis method finds all complex solutions
of a given system of zero dimensional polynomial equations when they are relatively small
systems.
Before discussing the semi-definite programming relaxation method,
we will study our discrete steady cavity flow problem
by the rational univariate representation \cite{rouillier}, \cite{noro-yokoyama},
which is a variation of the Gr\"obner basis method.
Results will be used to tune parameters of the sparse SDP relaxation method.

The $5 \times 5$ mesh is solvable with this method,
however the $5 \times 6$ mesh case is not solvable in one hour 
by current major implementations (Groebner(Fgb) in Maple 11, nd\_gr\_trace and tolex\_gsl in Risa/Asir). 
The system for the $5 \times 5$ mesh case contains $18$ variables and $9$ in the $18$ appear as linear and the other $9$ as quadratic variables. 

We sort the real solutions by a discretization of the kinetic energy
of the fluid ($\mbox{mass} \cdot \mbox{velocity}^2/2$), 
which is proportional to
\begin{equation}
 \int \int_{ABCD} \left( \frac{\partial \psi}{\partial y}\right)^2+\left(\frac{-\partial \psi}{\partial  x}\right) ^2 dx dy 
\label{energyintegral}
\end{equation}
By approximating integral (\ref{energyintegral}) via the central difference discretization,
we obtain our discrete energy function $F$, given as
\begin{equation}
\begin{array}{rl}
F(\psi,\omega)
& = \displaystyle\sum_{2\leq i,j\leq N-1} \left(\frac{(\psi_{i+1,j} - \psi_{i-1,j})^2}{4h^2}+ \frac{(\psi_{i,j+1} - \psi_{i,j-1})^2}{4h^2}\right)\; h^2 \\
& = \displaystyle\frac{1}{4}\sum_{2\leq i,j \leq N-1} \psi_{i+1,j}^2+\psi_{i-1,j}^2+\psi_{i,j+1}^2+\psi_{i,j-1}^2-2\psi_{i+1,j}\psi_{i-1,j}-2\psi_{i,j+1}\psi_{i,j-1}.
\end{array}
\label{objfun}
\end{equation}

\begin{figure}[ht]
\begin{center}
\includegraphics[width=0.32\textwidth, height=0.18\textheight]{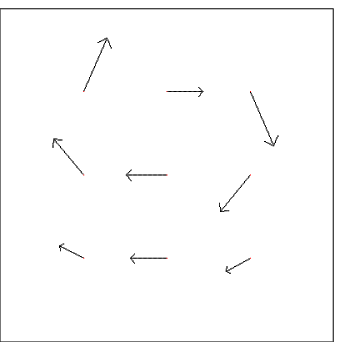} 
\includegraphics[width=0.32\textwidth, height=0.18\textheight]{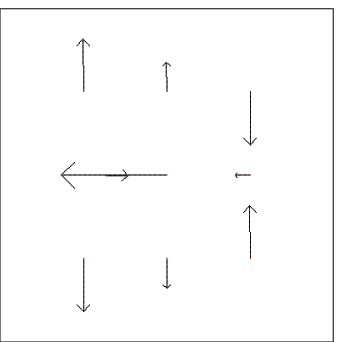} 
\includegraphics[width=0.32\textwidth, height=0.18\textheight]{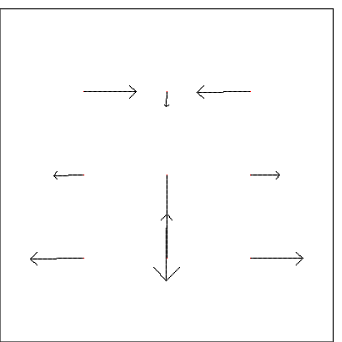} 
\caption{$R=0.01$, $v=1$, $N=5$, ${\rm M}=1,1/10^5, 1/10^5$.
There are $26$ real solutions.}
\label{r001n3v1}
\end{center}
\end{figure}

\begin{figure}[ht]
\begin{center}
\includegraphics[width=0.32\textwidth, height=0.18\textheight]{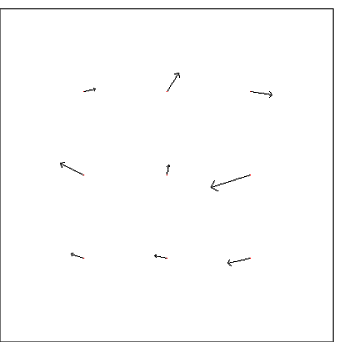} 
\includegraphics[width=0.32\textwidth, height=0.18\textheight]{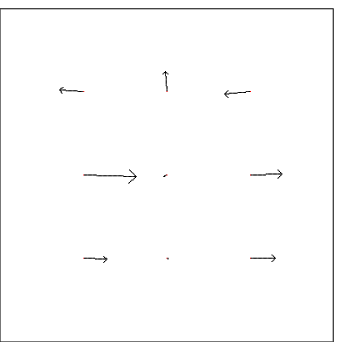} 
\includegraphics[width=0.32\textwidth, height=0.18\textheight]{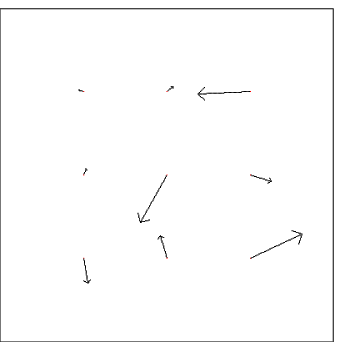} 
\caption{$R=500$, $v=1$, $N=5$, ${\rm M}=2,1/10, 1/10$.
There are $14$ real solutions.}
\label{r500n3v1}
\end{center}
\end{figure}

\begin{figure}[ht]
\begin{center}
\includegraphics[width=0.32\textwidth, height=0.18\textheight]{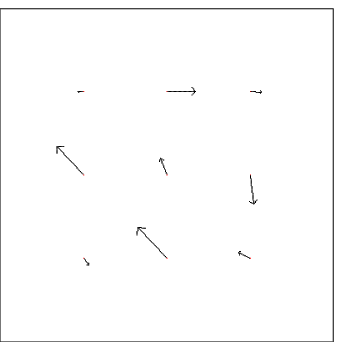} 
\includegraphics[width=0.32\textwidth, height=0.18\textheight]{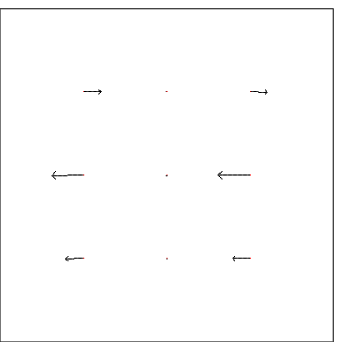} 
\includegraphics[width=0.32\textwidth, height=0.18\textheight]{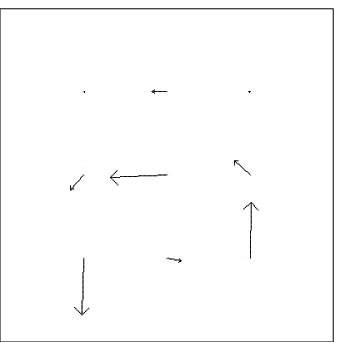} 
\caption{$R=40000$, $v=1$, $N=5$, ${\rm M}=4,2, 2$. 
There are $20$ real solutions.}
\label{r40000n3v1}
\end{center}
\end{figure}

\begin{figure}[ht]
\begin{center}
\includegraphics[width=0.32\textwidth, height=0.18\textheight]{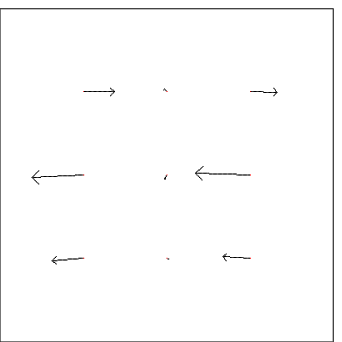} 
\includegraphics[width=0.32\textwidth, height=0.18\textheight]{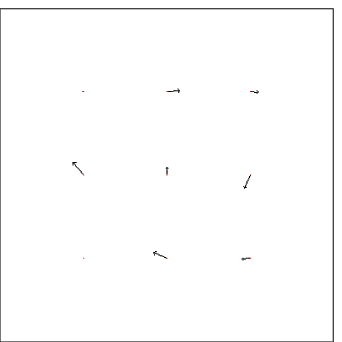} 
\includegraphics[width=0.32\textwidth, height=0.18\textheight]{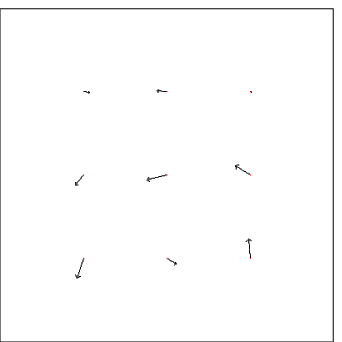} 
\caption{$R=500$, $v=10$, $N=5$, ${\rm Mag}=2,1/10, 1/10$.
There are $18$ real solutions.}
\label{r500n3v10}
\end{center}
\end{figure}

Figure \ref{r001n3v1} - \ref{r500n3v10} illustrate the approximations for some velocity vectors
$(\partial \psi/\partial y, - \partial \psi/\partial x)$
of solutions of the discrete steady cavity flow problem
obtained by the Gr\"obner basis method.
We sort real solutions by the discrete energy function.
The left solution is the minimum energy solution,
the center solution is the second energy solution,
the right solution is the 3rd energy solution.
$R$ is the Reynolds number, $v$ is the velocity of stream
along the boundary, $M$ is the magnification factor to display
the velocity vector, i.e. we display $(M \partial \psi/\partial y, -M \partial \psi/\partial x)$ in case $M\neq 1$.

It follows from these data that
the minimal energy and the energy gap between
the minimal energy and the second energy seem to decrease 
when $R$ increases.
It is interesting that the minimal energy solutions have a vortex
of clockwise direction, 
but some of 3rd energy solutions have a vortex of counterclockwise
direction, which are apparently fake solutions.
 
\section{Sparse SDP Relaxation Method}

The main contribution of this paper is to propose an algorithm that enumerates the smallest kinetic energy solutions of the discrete steady cavity flow problem  DSCF($R,v,N$) starting with the minimum energy solution. The key element of this algorithm is to apply the {\it sparse semidefinite program relaxation method} (SDPR) to solve the DSCF($R,v,N$). The SDPR for PDEs was proposed in \cite{mknt} and is based on the idea to take the polynomial system derived from a finite difference discretization of a differential equation and its boundary conditions (for instance: DSCF($R,v,N$)) as constraints for an optimization problem. After choosing a further polynomial function as the objective of the optimization problem, a polynomial optimization problem (POP) of the form
\begin{equation}
\begin{array}{lll}
\min & F(x)\\
\text{s.t. } & g_j(x) \geq 0 & \forall\; j\in\left\{1,\ldots,k\right\},\\
 & h_i(x) = 0 & \forall\; i\in\left\{ 1,\ldots,l \right\}.
\end{array}
\label{pop}
\end{equation}
is obtained. As shown in \cite{mknt}, polynomial optimization problems derived from differential equations satisfy structure sparsity patterns and the sparse SDP relaxations due to \cite{wkkm} can be applied to approximate the solution of POP (\ref{pop}). The crucial point is how to choose the objective function $F$ in POP (\ref{pop}). In case that several solutions to a discretized PDE problem exist, the choice of the objective function allows to select solutions of particular interest. For the cavity flow problem, we are interested in the solution which minimizes the kinetic energy (\ref{energyintegral}). Thus, for the cavity flow problem we yield the previously derived function $F$ (\ref{objfun}) as a canonical choice for the objective function of (\ref{pop}):
\begin{equation}
F(\psi,\omega)
 = \displaystyle\frac{1}{4}\sum_{2\leq i,j \leq N-1} \psi_{i+1,j}^2+\psi_{i-1,j}^2+\psi_{i,j+1}^2+\psi_{i,j-1}^2-2\psi_{i+1,j}\psi_{i-1,j}-2\psi_{i,j+1}\psi_{i,j-1}.
\end{equation}
We define the functions $g_{i,j}^1,\, g_{i,j}^2:\, \RR^{2N^2}\rightarrow\RR$ as
\begin{equation*}
\begin{array}{rl}
g_{i,j}^1(\psi,\omega) =& -4\omega_{i,j}+\omega_{i+1,j}+\omega_{i-1,j} + \omega_{i,j+1} + \omega_{i,j-1}\\
& + \frac{R}{4} (\psi_{i+1,j}-\psi_{i-1,j})(\omega_{i,j+1} - \omega_{i,j-1})\\
& - \frac{R}{4} (\psi_{i,j+1}-\psi_{i,j-1})(\omega_{i+1,j}-\omega_{i-1,j}),\\
g_{i,j}^2(\psi,\omega) = &-4\psi_{i,j}+\psi_{i+1,j}+\psi_{i-1,j}+\psi_{i,j+1}+\psi_{i,j-1}+h^2\omega_{i,j},
\end{array}
\end{equation*}
Taking into account DSCF($R,v,N$) and the objective function $F$, we derive the polynomial optimization problem,
\begin{equation}
\begin{array}{rll}
\min & F(\psi,\omega)\\
\text{s.t. } & g_{i,j}^1(\psi,\omega) = 0 & \forall\, 2\leq i,j\leq N-1,\\
& g_{i,j}^2(\psi,\omega) = 0 & \forall\, 2\leq i,j\leq N-1,\\
& \psi_{i,j} = 0 & \forall\, (i,j)\in\left\{1,N\right\}\times\left\{1,\ldots,N\right\}\cup\left\{1,\ldots,N\right\}\times\left\{1,N\right\}, \\
& \omega_{1,j} = -2\frac{\psi_{2,j}}{h^2} & \forall\, j\in \left\{ 1,\ldots,N \right\},\\
& \omega_{N,j} = -2\frac{\psi_{N-1,j}}{h^2} & \forall\, j\in \left\{ 1,\ldots,N \right\},\\ 
& \omega_{i,1} = -2\frac{\psi_{i,2}}{h^2} & \forall\, i\in \left\{ 1,\ldots,N \right\},\\
& \omega_{i,N} = -2\frac{\psi_{i,N-1} + v\,h}{h^2} & \forall\,i\in \left\{ 1,\ldots,N \right\}. 
\end{array}
\label{CFPOP}
\end{equation}
We call POP (\ref{CFPOP}) the {\bf steady cavity flow optimization problem} $\mathbf{CF(R,v,N)}$ with Reynold's number $R$, boundary velocity $v$ and discretization $N$ as parameters. As all polynomials in (\ref{CFPOP}) are of degree at most two,  $CF(R,v,N)$ is a {\bf quadratic optimization problem (QOP)}. In fact, a further classification is possible for $R=0$ and $R\neq 0$.

\begin{prop}
\begin{enumerate}
\item[a) ]
$CF(0,v,N)$ is a \textbf{convex quadratic program} for any $v$ and $N$.
\item[b) ]
$CF(R,v,N)$ is \textbf{non-convex} for any $v$ and $N$, if $R\neq 0$.
\end{enumerate}
\end{prop}
Proof:
\begin{enumerate}
\item[a) ]
In case $R=0$ all constraints are linear. Furthermore, the objective function can be written as $F=\sum_{i,j}F_{i,j}^1 + F_{i,j}^2$, where 
$$ F_{i,j}^1(\psi,\omega) = \left( \psi_{i-1,j},\; \psi_{i+1,j}\right) \;\left(\begin{array}{cc} 2 & -2 \\ -2 & 2\end{array}\right)\;\; \left( \begin{array}{c} \psi_{i-1,j} \\ \psi_{i+1,j}\end{array} \right).$$
It follows that $F_{i,j}$ is convex as $\left(\begin{array}{cc} 2 & -2 \\ -2 & 2\end{array}\right)$ positive semidefinite with eigenvalues $0$ and $4$. The convexity of $F_{i,j}^2$ follows analoguously. Thus, $F$ can be written as a sum of convex function and is therefore convex as well. The proposition follows.
\item[b) ]
In case $R\neq 0$, the equality constraint function $g_{i,j}^1$ is indefinite quadratic. Thus, $CF(R,v,N)$ is a non-convex quadratic program. $\qquad\square$
\end{enumerate}
It is our aim to apply the methods proposed in \cite{mknt} to solve $CF(R,v,N)$ and the underlying discrete steady cavity flow problem, i.e. to approximate the solutions of (\ref{CFPOP}) by solutions of a hierarchy of semidefinite program relaxations SDPR($w$) constructed in \cite{wkkm}, where $w$ denotes the order of the semidefinite program (SDP) relaxation. In theory, the solution of SDPR($w$) converges to the optimal solution for (\ref{CFPOP}) for $w\rightarrow\infty$. Nevertheless, the capacity of present SDP solvers restricts the choice of the relaxation order $w$, as the size of SDPR($w$) grows rapidly in $w$. However, as pointed out in \cite{wkkm} for many POPs it is suffcient to choose a relaxation order $w\in\left\{w_{\min},\ldots, w_{\min}+3\right\}$ to approximate the POP's minimizer accurately. For a general POP (\ref{pop}), $w_{\min}$ denotes the minimal relaxation order, which is given by
$$w_{\min} = \max\left\{\lceil \frac{\deg F}{2}\rceil ,\max_{1\leq j\leq k }\lceil \frac{\deg g_j}{2}\rceil,\max_{1\leq i\leq l}\lceil\frac{\deg h_i}{2}\rceil \right\}$$
Therefore, for $CF(R,v,N)$ holds $w_{\min}=1$.

\begin{rmk}
It is a well known result (c.f. \cite{lautut}), that SDPR(1) and (\ref{pop}) are equivalent, in case that the POP (\ref{pop}) is a convex quadratic program. Thus, solving $CF(0,v,N)$ is equivalent to solving a SDP.
\label{r0rmk}
Moreover, it is easy to show that the contraints admit only one feasible point when $R=0$.
\end{rmk}

\subsection{Tightening the SDP relaxation and improving the accuracy}
As stated before, the solution of SDPR($w$) converges to the optimizer of the POP for $w\rightarrow\infty$. Nevertheless, as the dimension $n$ of $CF(R,v,N)$ is given by $n=2(N-2)^2$, choosing a relaxation order $w$ greater than 2 for a medium scale discretization $N$ yields a SDP which requires too much memory in order to be solved by the used SDP-solver SeDuMi \cite{sedumi}. Therefore, we have to restrict ourselves to $w=1,2$ for small scale $N$, or even to $w=1$ for medium scale $N$. We cannot expect that SDPR(1) or SDPR(2) provide accurate approximations to the optimal solution for any $R$ and $v$. In order to tighten the SDP relaxation SDPR(1) and SDPR(2), respectively, we impose lower and upper bounds $\text{lbd}^{\psi}$, $\text{lbd}^{\omega}, \; \text{ubd}^{\psi}$ and $\text{ubd}^{\omega}\in\RR^{N^2}$ such that
\begin{equation}
\text{lbd}_i^\psi \leq \psi_i \leq \text{ubd}_i^\psi \quad\text{ and }\quad \text{lbd}_i^\omega \leq \omega_i \leq \text{ubd}_i^\omega\quad\forall\; 1\leq i\leq N^2
\label{lbdubd}
\end{equation}
holds.\\
Furthermore, we may apply additional locally convergent optimization techniques. For instance Newton's method for nonlinear systems can be applied to DSCF($R,v,N$) where the SDPR($w$) solution is taken as the starting point. Or alternatively, (\ref{CFPOP}) is approximated by sequential quadratic programming (SQP) \cite{boggs}, again, with the SDPR($w$) solution as starting point of the algorithm. Combining the sparse SDP relaxation with Newton's method or SQP is summarized in the scheme:

\begin{method} {\bf SDPR method}
\begin{enumerate}
\item 
Choose a boundary velocity $v$, grid discretization $N$ and Reynolds number $R'$.
\item
Apply SDPR($w$) to $CF(R',v,N)$ and obtain solution $\tilde{u}$.
\item
Apply sequential quadratic programming (SQP) to $CF(R',v,N)$ starting from $\tilde{u}$, and obtain $\hat{u}$.
\item
Apply Newton's method to DSCF($R',v,N$) starting from $\hat{u}$ or $\tilde{u}$, and obtain $u$.
\end{enumerate}
We note that  Step 3 and Step 4 are optional.
\label{sdprm}
\end{method}

\subsection{Enumeration algorithm for finding the $k$ smallest energy solutions}
As mentioned in section 2, we conjecture the number of solutions of the discrete steady cavity flow problem DSCF($R,v,N$) is finite, i.e. the feasible set of $CF(R,v,N)$ is finite. Method \ref{sdprm} enables us to approximate the global minimal solution $u^{\star} = u^{(1)\star}:=(\psi^{(1)\star},\omega^{(1)\star})$ of $CF(R,v,N)$. Beside the minimum solution, we are also interested in finding the solution $u^{(2)\star}$ with the second smallest kinetic energy, the solution $u^{(3)\star}$ with the third smallest kinetic energy or in general the solution $u^{(k)\star}$ with the $k$th smallest kinetic energy. Based on the SDPR method we propose an algorithm that enumerates the $k$ smallest kinetic energy solutions of $CF(R,v,N)$. Our algorithm shares the idea of separating the feasible set by additional constraints with Branch-and-Bound and cutting plane methods that are used for solving mixed integer linear programs and general concave optimization problems \cite{HorstPardalosThoai}. In contrast to the linear constraints of those methods we impose quadratic constraints to separate the feasible set. Moreover, $CF(R,v,N)$ is a non-convex continuous quadratic optimization problem for $R\neq 0$. It may be worth investigating in future how extensions of Branch-and-Cut methods for certain nonconvex problems \cite{fukuda} can be used in our setting.

\begin{alg}
Enumerate the $k$ smallest solutions:\\
Given $u^{(k-1)}$, the approximation to  the $(k-1)$th energy solution obtained by solving $\text{SDPR}^{k-1}(w)$.
\begin{enumerate}
\item
Choose $\epsilon_1^k$ and $\epsilon_2^k>0$.
\item
Choose integers $b_1^k,\; b_2^k\in\left\{ 1,\ldots,(N-2)^2 \right\}$. 
\item
Add the following quadratic constraints to $\text{SDPR}^{k-1}(w)$ and denote the resulting (tighter) SDP relaxation as $\text{SDPR}^k(w)$.
\begin{equation}
\begin{array}{ll}
(u_j - u_j^{k-1})^2 := (\psi_j - \psi_j^{k-1})^2 \geq \epsilon_1^k & \forall 1\leq j\leq b_1^k,\\
(u_{j+(N-2)^2} - u_{j+(N-2)^2}^{k-1})^2 := (\omega_j - \omega_j^{k-1})^2 \geq\epsilon_2^k & \forall 1\leq j\leq b_2^k.
\end{array}
\label{addNormCon}
\end{equation} 
\item
Solve $\text{SDPR}^k(w)$ with $w=1,2$ or larger, if possible. Obtain first approximation $u^{\text{SDP}(k)}$.
\item
Apply a local optimization technique as for instance Newton's method or SQP with $u^{\text{SDP}(k)}$ as starting point. Obtain $u^{(k)}$ as an approximation to $u^{(k)\star}$.
\item
Iterate steps 1--5.
\end{enumerate}
\label{enumSolAlg}
\end{alg}
The idea of Algorithm \ref{enumSolAlg} is to impose an additional polynomial inequality constraint 
\begin{equation*}
\begin{array}{ll}
(u_j - u_j^{k-1})^2 := (\psi_j - \psi_j^{k-1})^2 \geq \epsilon_1^k & \forall 1\leq j\leq b_1^k,\\
(u_{j+(N-2)^2} - u_{j+(N-2)^2}^{k-1})^2 := (\omega_j - \omega_j^{k-1})^2 \geq\epsilon_2^k & \forall 1\leq j\leq b_2^k.
\end{array}
\end{equation*}
to the POP (\ref{CFPOP}) in iteration $k$, that excludes the solution $u^{k-1}$ from the feasible set of POP (\ref{CFPOP}) which was obtained in the previous iteration. In case that the feasible set of (\ref{CFPOP}) is finite and $u^{k-1}$ is sufficiently close to $u^{(k-1)\star}$, the new constraint excludes $u^{(k-1)\star}$ from the feasible set of (\ref{CFPOP}) and $u^{(k)\star}$ is the new global minimizer of (\ref{CFPOP}). Of course, there are various alternatives to step 3 in Algorithm \ref{enumSolAlg}, in order to exclude $u^{(k-1)\star}$ from the POP's feasible set. One alternative constraint is
\begin{equation}
\left( u_i - u_i^{(k-1)\star}\right)u_{n+i} - \epsilon_i = 0\quad (1\leq i\leq b),
\label{altCon}
\end{equation}
where $b\in\left\{ 1,\ldots,n\right\},\;\epsilon_i>0$ and $u_{n+i}$ a new additional variable bounded by $-1$ and $1$. It is easy to see that (\ref{altCon}) is violated, if $u=u^{(k-1)\star}$. Nevertheless, it turned out that the numerical performance of (\ref{altCon}) is inferior to the one of (\ref{addNormCon}) for our problem DSCF($R,v,N$). The right tuning of parameters $\epsilon_i$ and $b$ is far more difficult for (\ref{altCon}) compared to (\ref{addNormCon}). A second alternative to exclude $u^{(k-1)\star}$ are $l_p$-norm constraints such as
\begin{equation}
\parallel u-u^{\star} \parallel_p = \left(\sum_{i=1}^n \left( u_i-u_i^{(k-1)} \right)^p\right)^{\frac{1}{p}}\geq\epsilon,
\label{lpnormCon}
\end{equation}
for $p\geq 1$. The disadvantage of the constraints (\ref{lpnormCon}) is, they destroy the sparsity of the POP (\ref{CFPOP}), as all $u_i\; (i=1,\ldots,n)$ occur in the same constraint. Therefore the advantage of the sparse SDP relaxations is lost and the POP can not be solved efficiently anymore. An exception is the infinity norm constraint
\begin{equation}
\parallel u-u^{\star} \parallel_{\infty} = \max_{1\leq i\leq n}\mid u_i-u_i^{(k-1)} \mid\geq\epsilon.
\label{infnormCon}
\end{equation}
In fact, the infinity norm constraints (\ref{infnormCon}) are equivalent to the proposed constraints (\ref{addNormCon}) for $b_1=b_2=\frac{n}{2}$. Because it preserves sparsity and its numerical performance is better than the one of (\ref{altCon}), we impose (\ref{addNormCon}) as additional constraints in Algorithm \ref{enumSolAlg}.  We obtain the following results for Algorithm \ref{enumSolAlg}.
\begin{prop}
Let $R,\; v$ and $N$ be fixed. Let $(u^{(1)},\ldots,u^{(k-1)})$ be the output of the first $(k-1)$ iterations of Algorithm (\ref{enumSolAlg}). If this output is a sufficiently close approximation of the vector of $(k-1)$ smallest kinetic energy solutions $(u^{(1)\star},\ldots,u^{(k-1)\star})$, and if the feasible set of POP (\ref{CFPOP}) is finite, then there exist $b\in\left\{1,\ldots,n\right\}$ and $\epsilon\in\RR^b$  such that  the output $u^{(k)}$ of Algorithm \ref{enumSolAlg}  (for $k$th iteration) satisfies
\begin{equation}
u^{(k)}(w) \rightarrow u^{(k)\star} \quad \text{ when } w\rightarrow\infty.
\end{equation}
\label{convThm}
\end{prop}
\paragraph{Proof: }
As each $u^{(j)}$ is in a neighborhood of $u^{(j)\star}$ for all $j\in\left\{1,\ldots,k-1\right\}$, we can choose $b\in\left\{ 1,\ldots,n \right\}$ and a vector $\epsilon\in\RR^b$, s.t. 
$$\forall j\in\left\{ 1,\ldots, k-1\right\}\;\exists i\leq b\text{ s.t. } \left( u_i-u_i^{(j)\star} \right)^2 < \epsilon_i,$$
and $$ \left( u_i-u_i^{(l)\star} \right)^2\geq\epsilon_i \;\forall l\geq k\;\forall i\in\left\{ 1,\ldots, b\right\}.$$
Let $CF(R,v,N)^{(k)}$ denote the $CF(R,v,N)$ with the $k$ systems of additional constraints given by step 3 in Algorithm \ref{enumSolAlg}, where the $k$th constraints are given by (\ref{addNormCon}) for the constructed $b$ and $\epsilon$. Then it holds
$$\text{feas}\left(CF(R,v,N)^{(k)}\right)=\text{feas}\left(CF(R,v,N)\right)\setminus\left\{ u^{(1)\star},\ldots,u^{(k-1)\star}\right\}.$$
Thus, $u^{(k)\star}$ is the global minimizer of $CF(R,v,N)^{(k)}$ and the global minimum is $F(u^{(k)\star})$. As the bounds (\ref{lbdubd}) guarantee the compactness of the feasible set, it holds with the convergence theorem for the sparse SDP relaxations \cite{lass2}
\begin{equation}
\begin{array}{lcll}
\min \text{SDPR}^{(k)}(w) & \rightarrow & \min CF(R,v,N)^{(k)}=F(u^{(k)\star}) & \text{ for } w\rightarrow\infty,\\
u^{(k)}(w) & \rightarrow & u^{(k)\star} & \text{ for } w\rightarrow\infty,
\end{array}
\label{convThmconv}
\end{equation}
under the assumption that $F(u^{(1)\star})<F(u^{(2)\star})<\ldots$, i.e. $u^{(k)\star}$ is the unique minimizer of $CF(R,v,N)^{(k)}.\qquad\square$

Although we have proven convergence, the capacity of current SDP solvers restricts the choice of the relaxation order $w$ to small integers, in our application typically to $w=1$ or $w=2$. Furthermore, we need to choose the parameters $\epsilon$ and $b$ appropriately, to obtain good approximations of the $k$ smallest kinetic energy solutions. In the following numerical examples we will discuss this question and show heuristics how to tune these two parameters.

\subsection{Numerical results}
We will demonstrate the numerical performance of the SDPR($w$) and Algorithm \ref{enumSolAlg} to enumerate the $k$ smallest solutions. All calculations are conducted on a Mac OS X with CPU 2.5GHz and 2 GB Memory. As an implementation of the sparse SDP relaxation we use the software SparsePOP \cite{sparsepop} and MATLAB optimization toolbox for standard SQP routines in order to improve the accuracy of the solution provided by SDPR($w$). The total accumulated computation time in seconds is denoted by $t_C$, the scaled feasibility error of a SDPR solution $u'$ w.r.t. the constraints of $CF(R,v,N)$ by $\epsilon_{\text{sc}}$.

\subsubsection{CF(4000,1,5)}

\begin{figure}[ht]
\begin{center}
\includegraphics[width=0.32\textwidth, height=0.18\textheight]{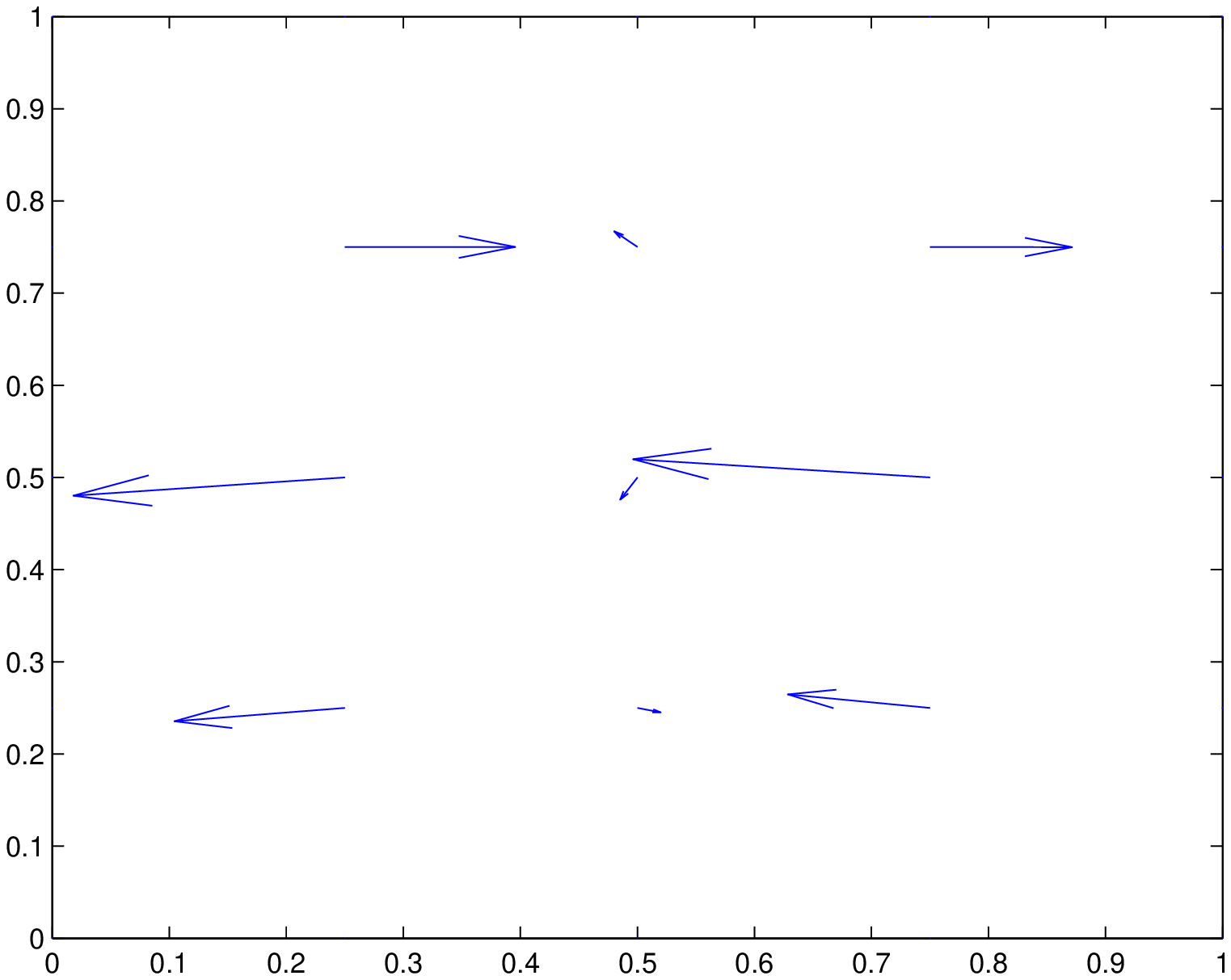} 
\includegraphics[width=0.32\textwidth, height=0.18\textheight]{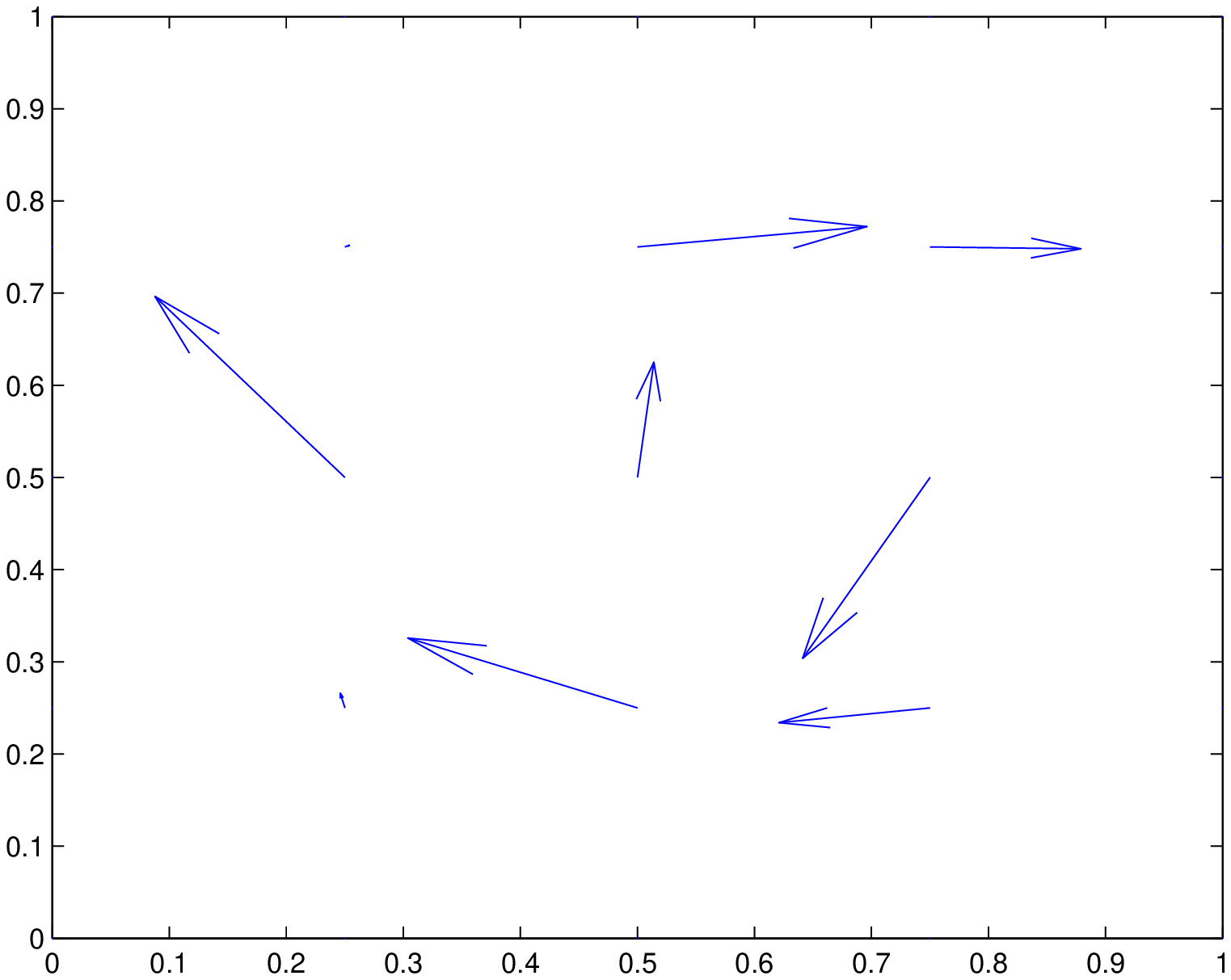} 
\includegraphics[width=0.32\textwidth, height=0.18\textheight]{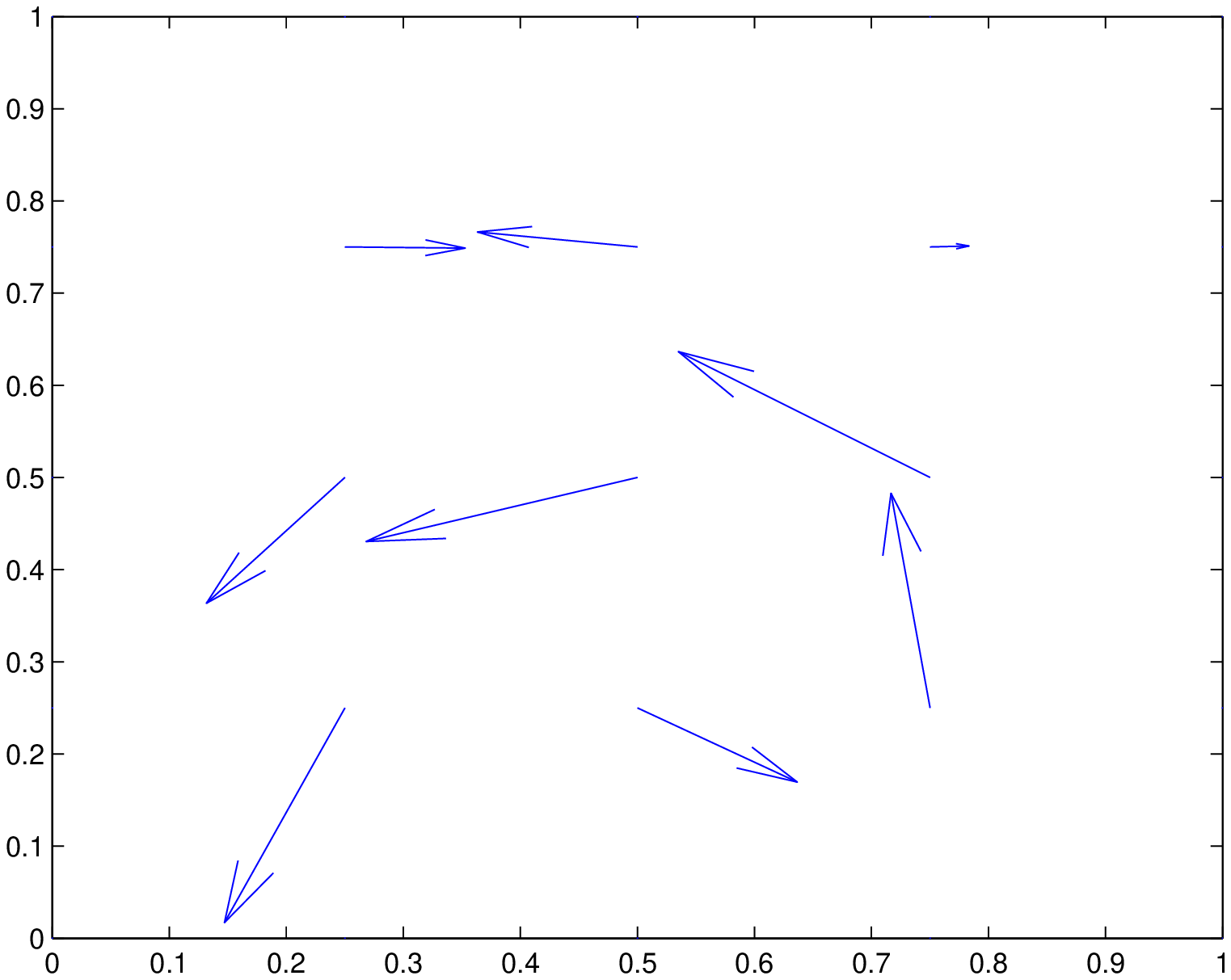} 
\caption{Interior of the solutions $u^{(0)}$ (left), $u^{(1)}$ (center) and $u^{(2)}$ (right) for $CF(4000,1,5)$}
\label{CF4000g5fig}
\end{center}
\end{figure}

In a first setting we choose the discretization $N=5$, i.e. the dimension of the POP (\ref{CFPOP}) is $n=2\cdot3^2 = 18$. This dimension is small enough to apply the polyhedral homotopy method \cite{hom4ps}
and the Gr\"obner basis method (Section \ref{section:gbmethod})  to determine all complex solutions of DSCF($R,v,N$). Therefore, we are able to verify whether the solutions provided by Algorithm \ref{enumSolAlg} are optimal. The computational results are given in Table \ref{CF4000g5res}.
\begin{table}[ht]
\begin{center}
\begin{tabular}{|r|r|r|r|r|r|r|r|r|r|}
\hline $k$ & $w$ & $\epsilon_1^k$ & $\epsilon_2^k$ & $b_1^k$ & $b_2^k$ & $t_C$ & $\epsilon_{\text{sc}}$ & $F(u^{(k)})$ & solution\\
\hline 0 & 1 & - & - & - & - & 2 & 2e-10 & 4.6e-4 & $u^{(0)}$\\
\hline 1 & 1 & 1e-3 & - & 3 & 0 & 5 & 5e-4 & 6.3e-4 & $u^{(1)}$\\
\hline 2 & 1 & 1e-3 & - & 3 & 0 & 8 & 5e-4 & 1.0e-3 & $u^{(2)}$\\
\hline
\end{tabular}
\caption{Results of Algorithm \ref{enumSolAlg} for $CF(4000,1,5)$}
\label{CF4000g5res}
\end{center}
\end{table}
Comparing our SDPR results to all solutions of the polynomial system obtained by polyhedral homotopy method or Gr\"obner basis method, it turns out that the solutions $u^{(0)}$, $u^{(1)}$ and $u^{(2)}$ indeed coincide with the three smallest energy solutions $u^{(0)\star}$, $u^{(1)\star}$ and $u^{(2)\star}$. All three solutions are pictured in Figure \ref{CF4000g5fig}. Note, that the third smallest energy solution $u^{(2)}$ shows a vortex in counter-clockwise direction, which may indicate that this solution is a fake solution.

\subsubsection{CF(20000,1,7)}

\begin{figure}[ht]
\begin{center}
\includegraphics[width=0.32\textwidth, height=0.22\textheight]{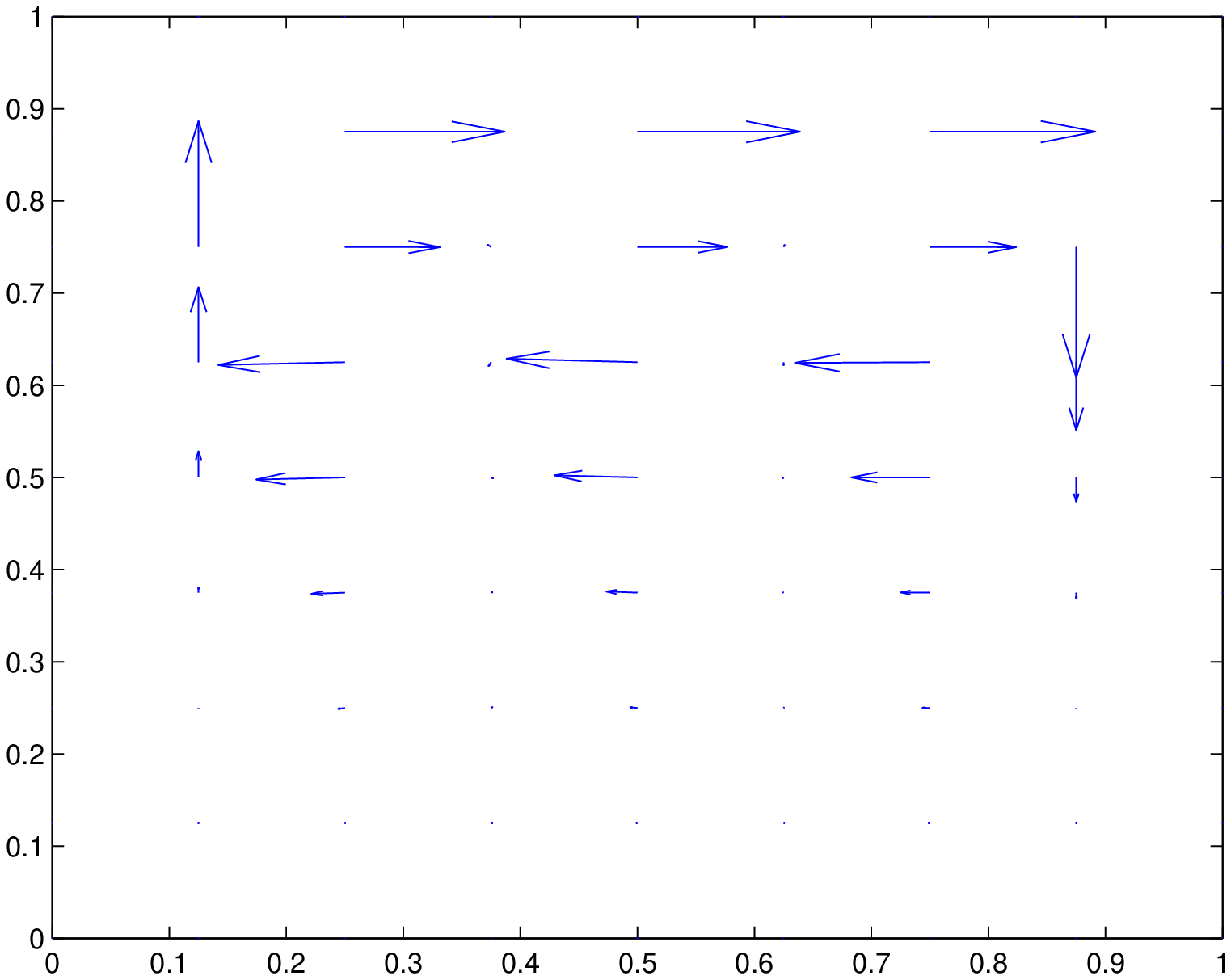} 
\includegraphics[width=0.32\textwidth, height=0.22\textheight]{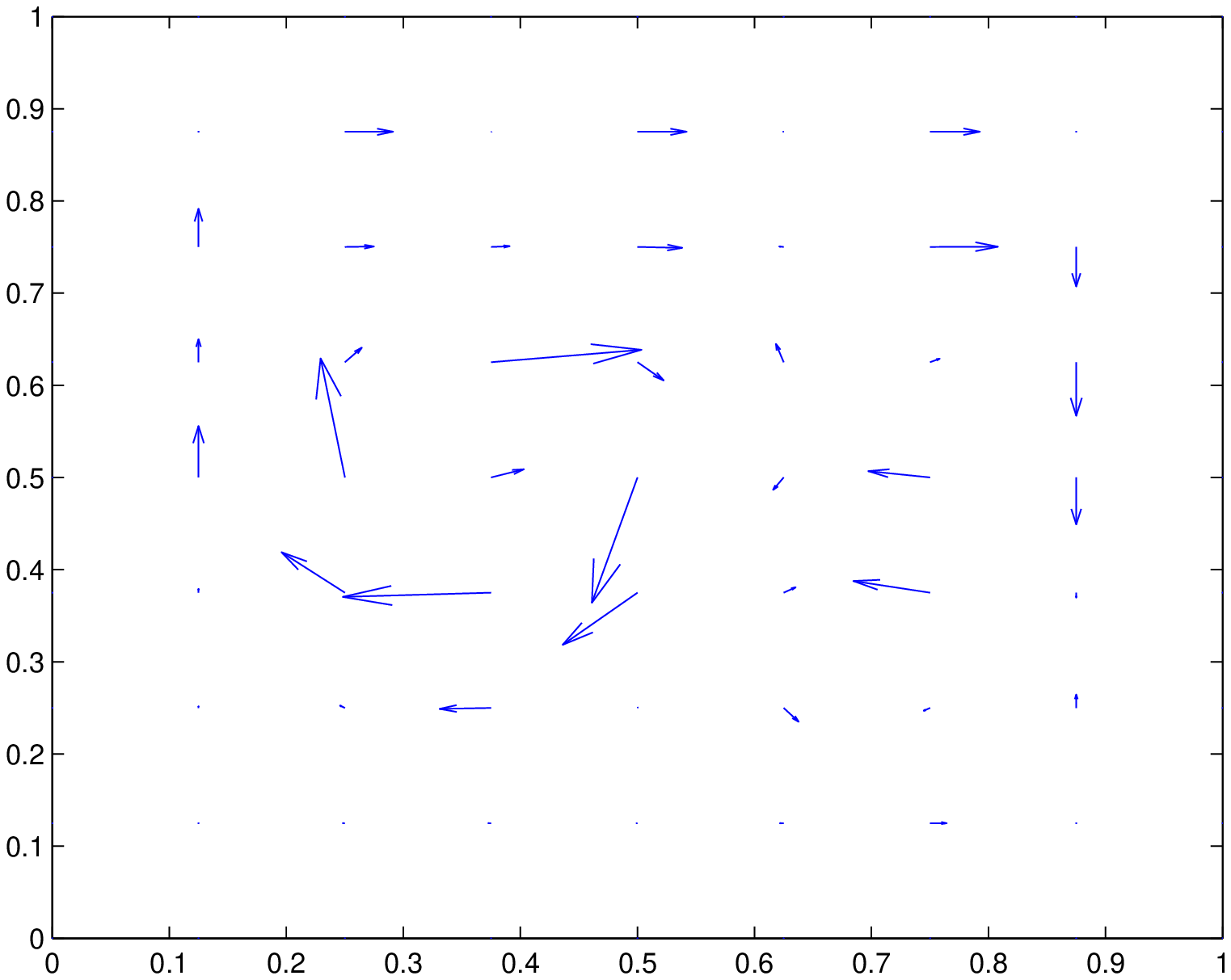} 
\includegraphics[width=0.32\textwidth, height=0.22\textheight]{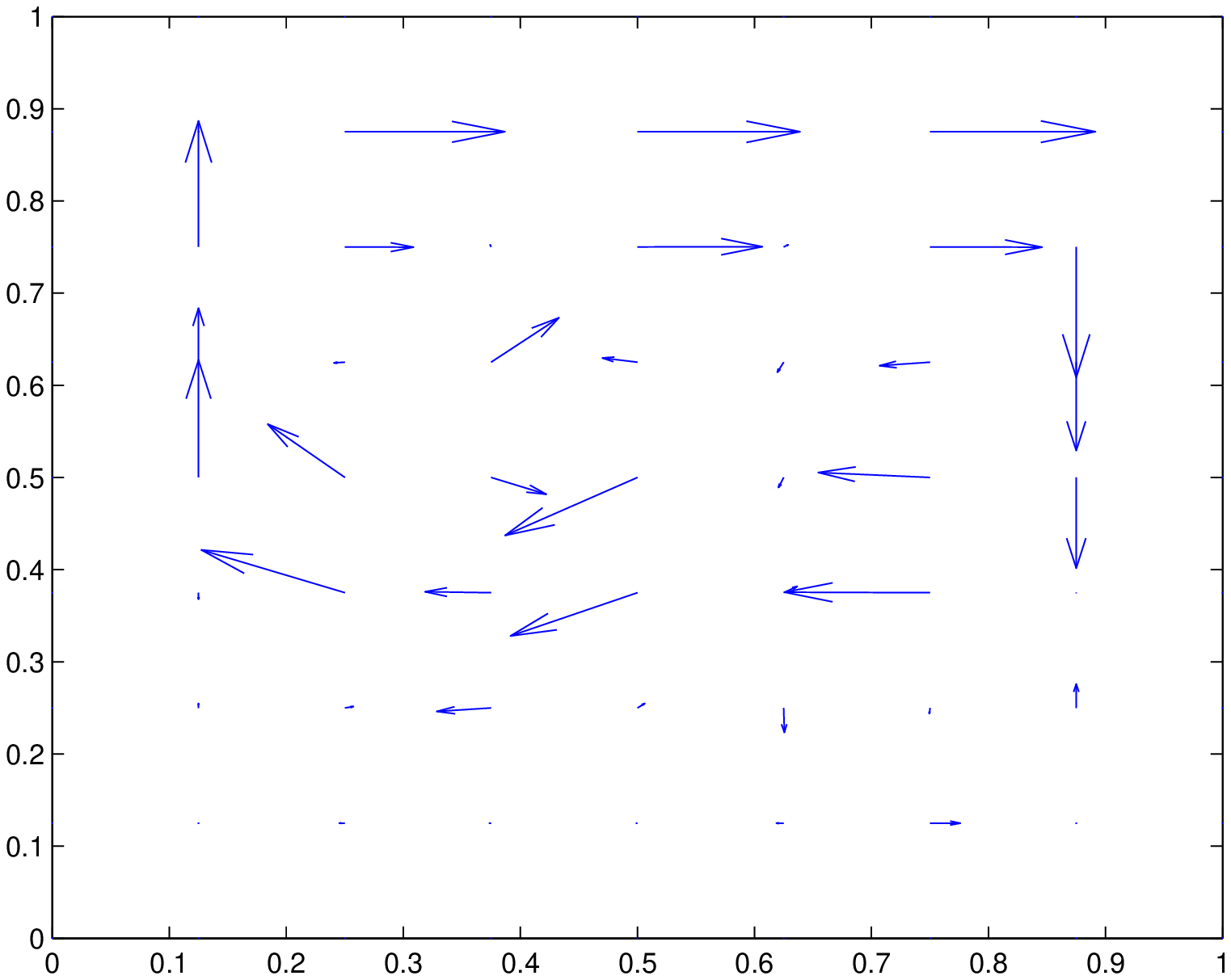} 
\caption{Solutions $u^{(0)}$ (left), $u^{(1)}$ (center) and $u^{(2)}$ (right) for $CF(20000,1,7)$}
\label{CF20000g7fig}
\end{center}
\end{figure}

We restrict ourselves to SDPR(1) for solving $CF(20000,1,7)$, as the size of the SDP relaxation with order $w=2$ resulting from this POP of dimension $n = 2\cdot 25^2=50$ is already too large to be solved in reasonable time. The computational results for applying Algorithm \ref{enumSolAlg} for different choices of the algorithm parameters are enlisted in Table \ref{CF20000g7anres}. The two parameter settings $(\epsilon_1^1,b_1^1)=(1e-3,\,1)$ and $(\epsilon_1^1,b_1^1)=(1e-6,\,5)$ are not sufficient to obtain an other solution than $u^{(0)}$, whereas $(\epsilon_1^1,b_1^1)=(1e-5,\,5)$ yields $u^{(1)}$, a solution of larger energy. After another iteration with $(\epsilon_1^2,b_1^2)=(1e-5,\,5)$ we obtain a third solution $u^{(3)}$ of even larger energy. The three solutions are pictured in Figure \ref{CF20000g7fig}.

\begin{table}[ht]
\begin{center}
\begin{tabular}{|r|r|r|r|r|r|r|r|r|r|}
\hline $k$ & $w$ & $\epsilon_1^k$ & $\epsilon_2^k$ & $b_1^k$ & $b_2^k$ & $t_C$ & $\epsilon_{\text{sc}}$ & $F(u^{(k)})$ & solution\\
\hline 0 & 1 & - & - & - & - & 2 & 3e-7 & 3.4e-4 & $u^{(0)}$\\
\hline 1 & 1 & 1e-3 & - & 1 & 0 & 5 & 5e-4 & 3.4e-4 & $u^{(0)}$\\
\hline 1 & 1 & 1e-6 & - & 5 & 0 & 5 & 6e-6 & 3.4e-4 & $u^{(0)}$\\
\hline 1 & 1 & 1e-5 & - & 5 & 0 & 9 & 5e-6 & 5.9e-4 & $u^{(1)}$\\
\hline 2 & 1 & 1e-5 & - & 5 & 0 & 14 & 5e-6 & 5.2e-3 & $u^{(2)}$\\
\hline
\end{tabular}
\caption{Results of Algorithm \ref{enumSolAlg} for $CF(20000,1,7)$}
\label{CF20000g7anres}
\end{center}
\end{table}

It is interesting to observe that $u^{(1)}$ and $u^{(2)}$ are one-vortex solutions, whereas there seems to be no vortex in the smallest energy solution $u^{(0)}$.

\subsubsection{CF(40000,1,7)}

Next, we examine $CF(40000,1,7)$, which is a good example to demonstrate that solving DSCF($R,v,N$) and POP (\ref{CFPOP}) is becoming more difficile for larger Reynolds numbers. As for the previous problem, the dimension of the POP is $n=50$, which is too large to be solved by Gr\"obner basis or polyhedral homotopy method. Our computational results are reported in Table \ref{CF40000g7anres}.

\begin{table}[ht]
\begin{center}
\begin{tabular}{|r|r|r|r|r|r|r|r|r|r|}
\hline $k$ & $w$ & $\epsilon_1^k$ & $\epsilon_2^k$ & $b_1^k$ & $b_2^k$ & $t_C$ & $\epsilon_{\text{sc}}$ & $F(u^{(k)})$ & solution\\
\hline 0 & 1 & - & - & - & - & 3 & 2e-7 & 3.4e-4 & $u^{(0)}(1)$ \\
\hline 1 & 1 & 5e-6 & - & 5 & 0 & 7 & 6e-9 & 7.3e-4 & $u^{(1)}(1)$ \\
\hline 2 & 1 & 5e-6 & - & 5 & 0 & 11 & 3e-6 & 5.9e-4 & $u^{(2)}(1)$ \\
\hline 3 & 1 & 8e-6 & - & 5 & 0 & 16 & 5e-6 & 2.3e-4 & $u^{(3)}(1)$ \\
\hline 0 & 2 & - & - & - & - & 5872 & 8e-10 & 2.6e-4 & $u^{(0)}(2)$ \\
\hline
\end{tabular}
\caption{Results of Algorithm \ref{enumSolAlg} for $CF(40000,1,7)$}
\label{CF40000g7anres}
\end{center}
\end{table}

\begin{figure}[ht]
\begin{center}
\includegraphics[width=0.32\textwidth, height=0.22\textheight]{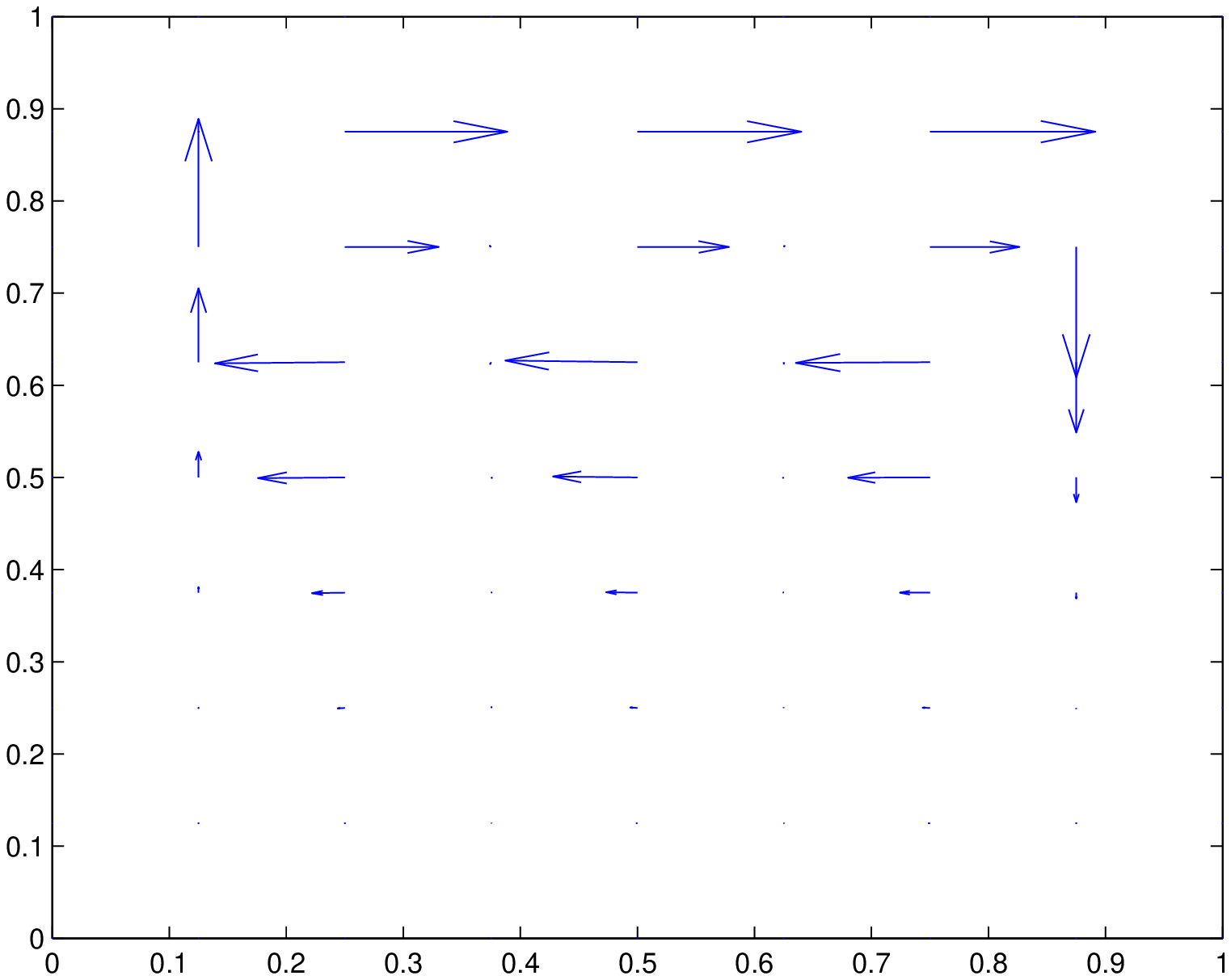} 
\includegraphics[width=0.32\textwidth, height=0.22\textheight]{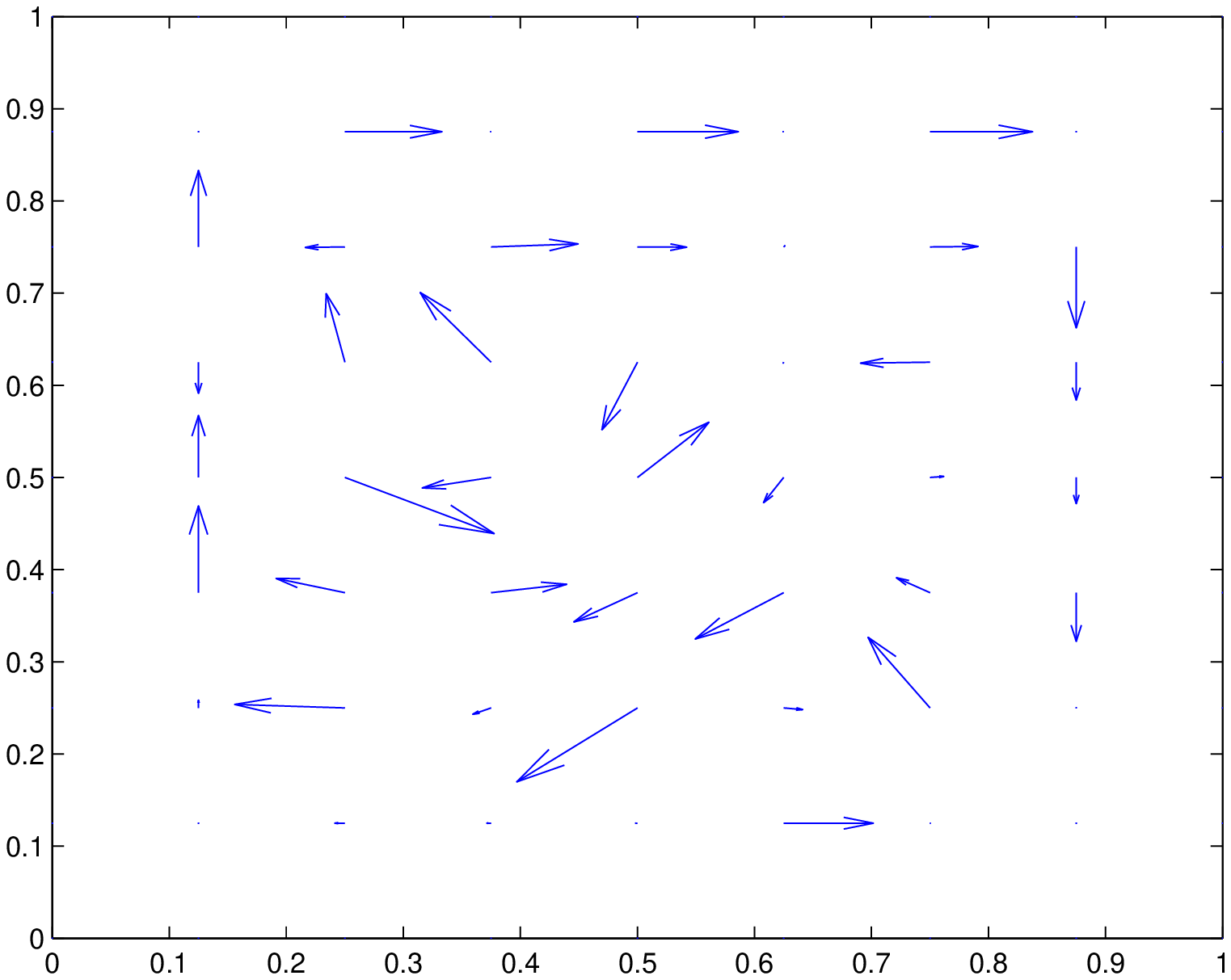} 
\caption{Solutions $u^{(0)}(1)$ (left) and $u^{(0)}(2)$ (right)}
\label{CF40000g7fig}
\end{center}
\end{figure}

\begin{figure}[ht]
\begin{center}
\includegraphics[width=0.32\textwidth, height=0.22\textheight]{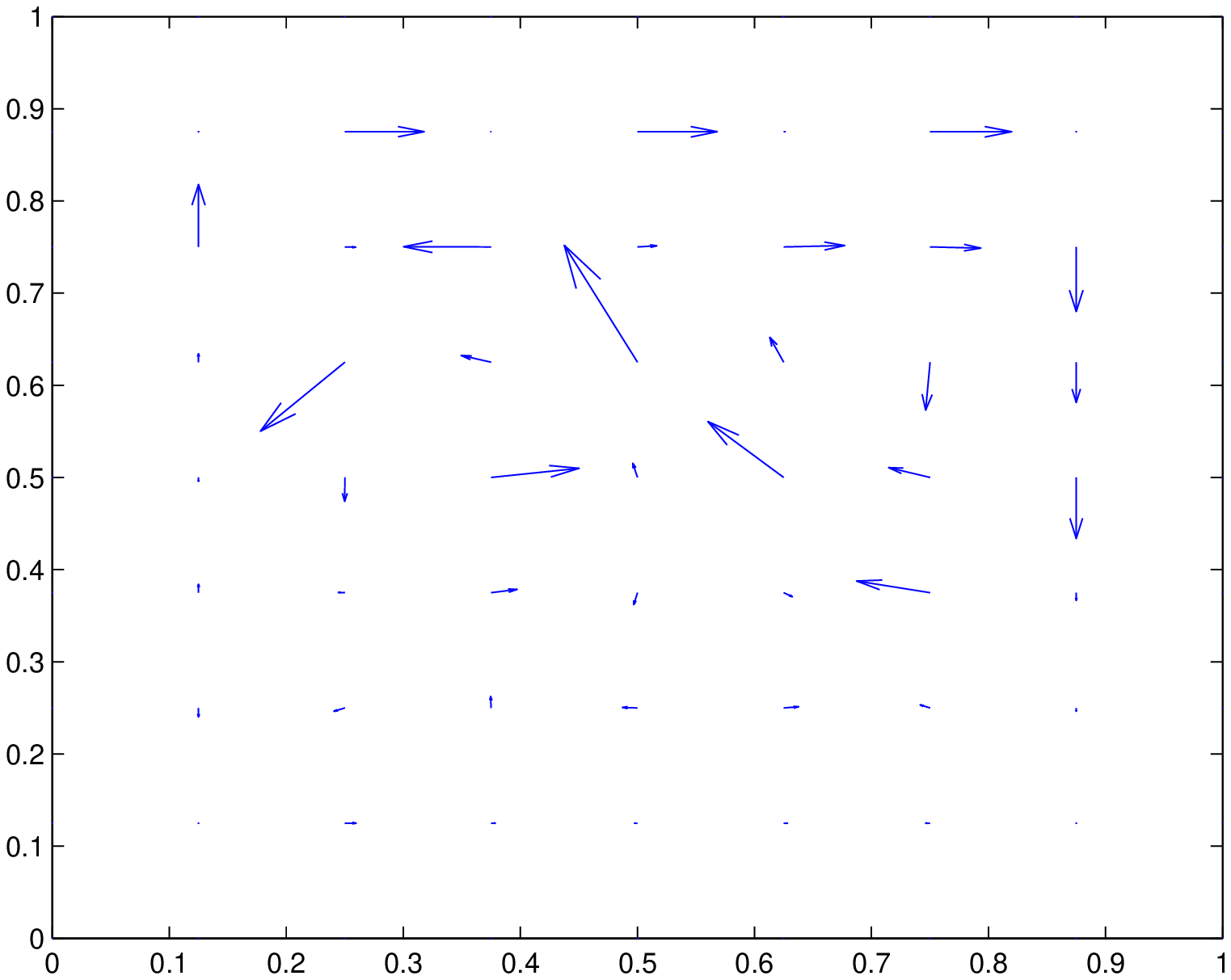} 
\includegraphics[width=0.32\textwidth, height=0.22\textheight]{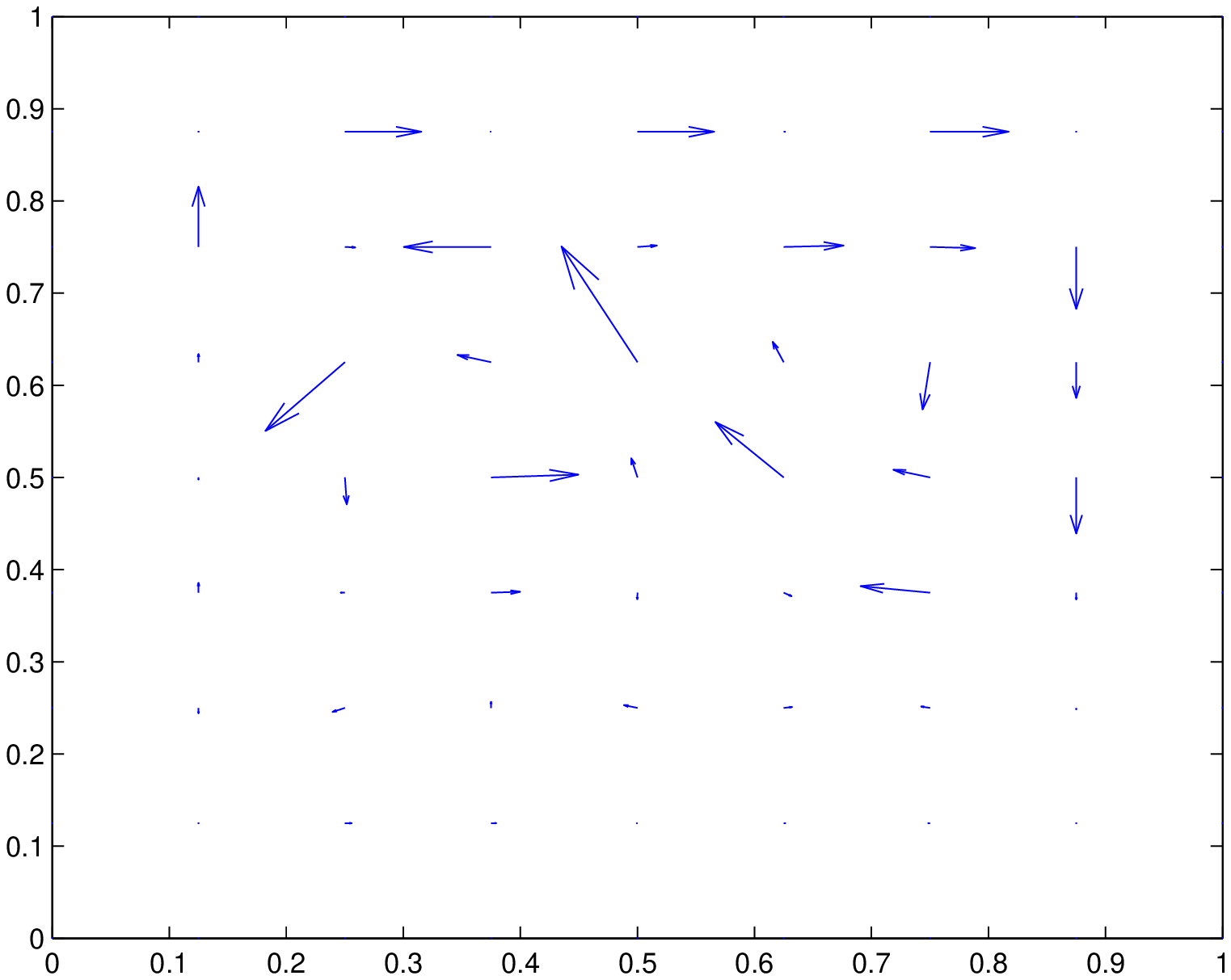} 
\includegraphics[width=0.32\textwidth, height=0.22\textheight]{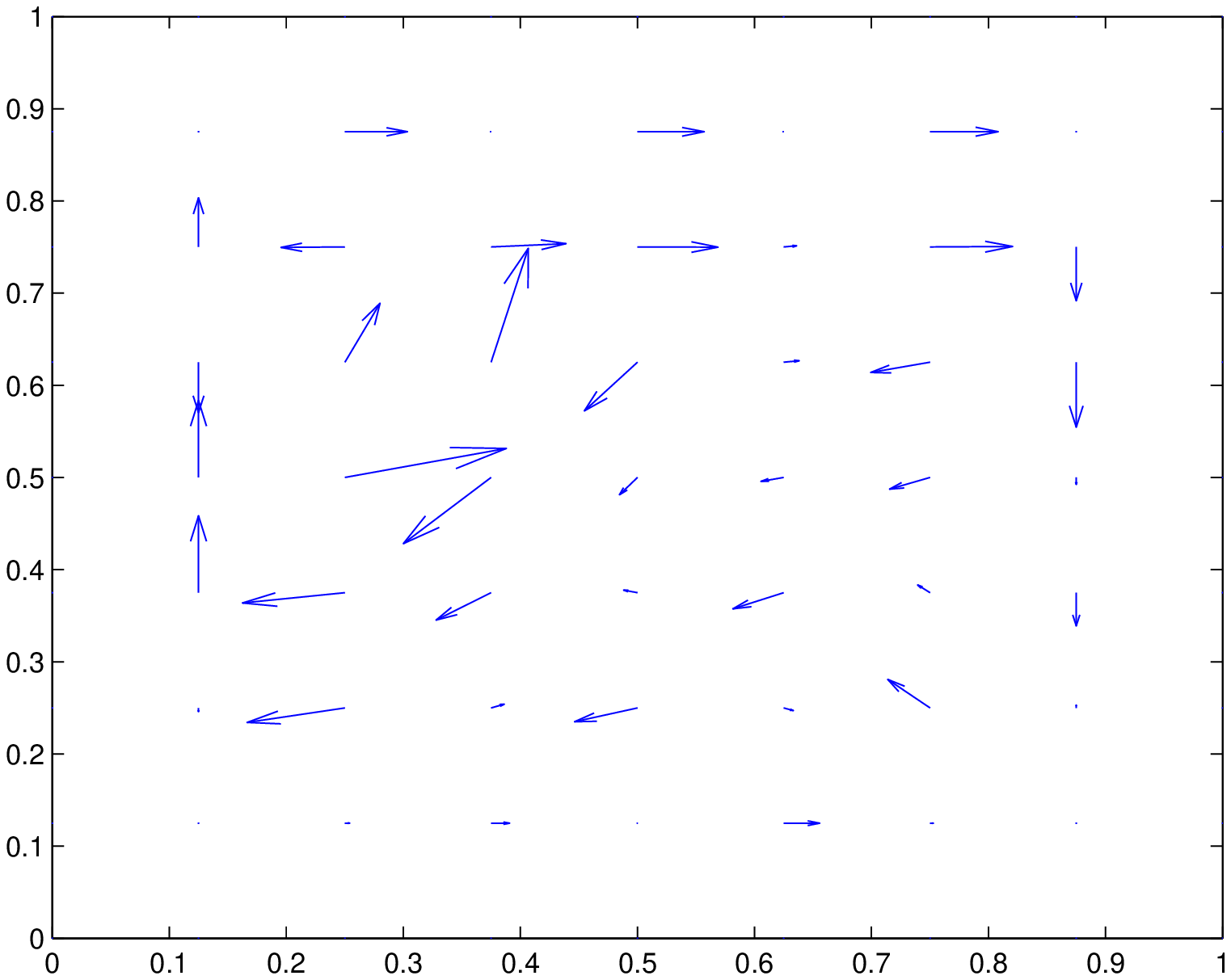} 
\caption{Solutions $u^{(1)}(1)$ (left), $u^{(2)}(1)$ (center) and $u^{(3)}(1)$ (right)}
\label{CF40000g7fig2}
\end{center}
\end{figure}

Solution $u^{(2)}(1)$ is of smaller energy than $u^{(1)}(1)$, and $u^{(3)}(1)$ is even of smaller energy than $u^{(0)}(1)$. This phenomenon can be explained by the fact, that the SDP relaxation with $w=1$ is not tight enough to yield a solution that converges to $u^{\star}$ under the local optimization procedure. The energy of $u^{(0)}(2)$ obtained by SDPR(2) is smaller than the one of $u^{(0)}(1)$, but it is not the global minimizer as well. In fact, Algorithm \ref{enumSolAlg} with SDPR(1) generates a better solution $u^{(3)}(1)$ (with smaller energy) in 3 iterations requiring 16 seconds computation time, compared to solution $u^{(0)}(2)$ obtained by applying SDPR(2)  to $CF(40000,1,7)$ requiring 5872 seconds. Thus, applying our enumeration algorithm with relaxation order $w=1$ is far more efficient than the original sparse SDP relaxation with $w=2$ for approximating the global minimizer of POP (\ref{CFPOP}) in this example; our enumeration algorithm fails with the relaxization order $w=1$, but we obtain the global minimizer efficiently by accident. The 5 different solutions are illustrated in Figure \ref{CF40000g7fig} and \ref{CF40000g7fig2}.

\section{Convergence of solutions in case of large Reynolds numbers $R$}

\subsection{SDPR($w$) for increasing discretization $N$}

In our previous experiments we derived small or even minimal energy solutions by Method \ref{sdprm} and Algorithm \ref{enumSolAlg} for various choices of the problem parameters $R$ and $v$ with discretization $N\in\left\{ 5,\, 6,\, 7 \right\}$. In case that we succeed, applying SDPR($w$) to $CF(R,v,N)$ yields the minimum kinetic energy solution $u^{\star}$ of the discrete steady cavity flow problem. The important question arises whether it is possible to expand these coarse grid minimum kinetic energy solutions $u^{\star}$ to finer grids with larger discretization $N$, i.e. whether these coarse grid solutions converge to analytic solutions of the original (continuous) steady cavity flow problem for $N\rightarrow\infty$. As pointed out in, e.g., \cite{cheng}, in case of  larger and larger Reynolds number $R$ and velocity $v$ the discrete steady cavity flow problem has to be solved for finer and finer grids, in order to obtain solutions converging to solutions of the continuous problem for $N\rightarrow\infty$. In this section we will adress the difficult problem to find solutions for $CF(R,v,N)$ converging to continuous solutions for large $R$ and pose a related question and a conjecture based on computational experiments. As in section 4, the calculations are conducted on a Mac OS X with CPU 2.5GHz and 2 GB Memory. 

\subsubsection{$CF(100,1,N)$}
\begin{figure}[ht]
\begin{center}
\includegraphics[width=0.32\textwidth, height=0.16\textheight]{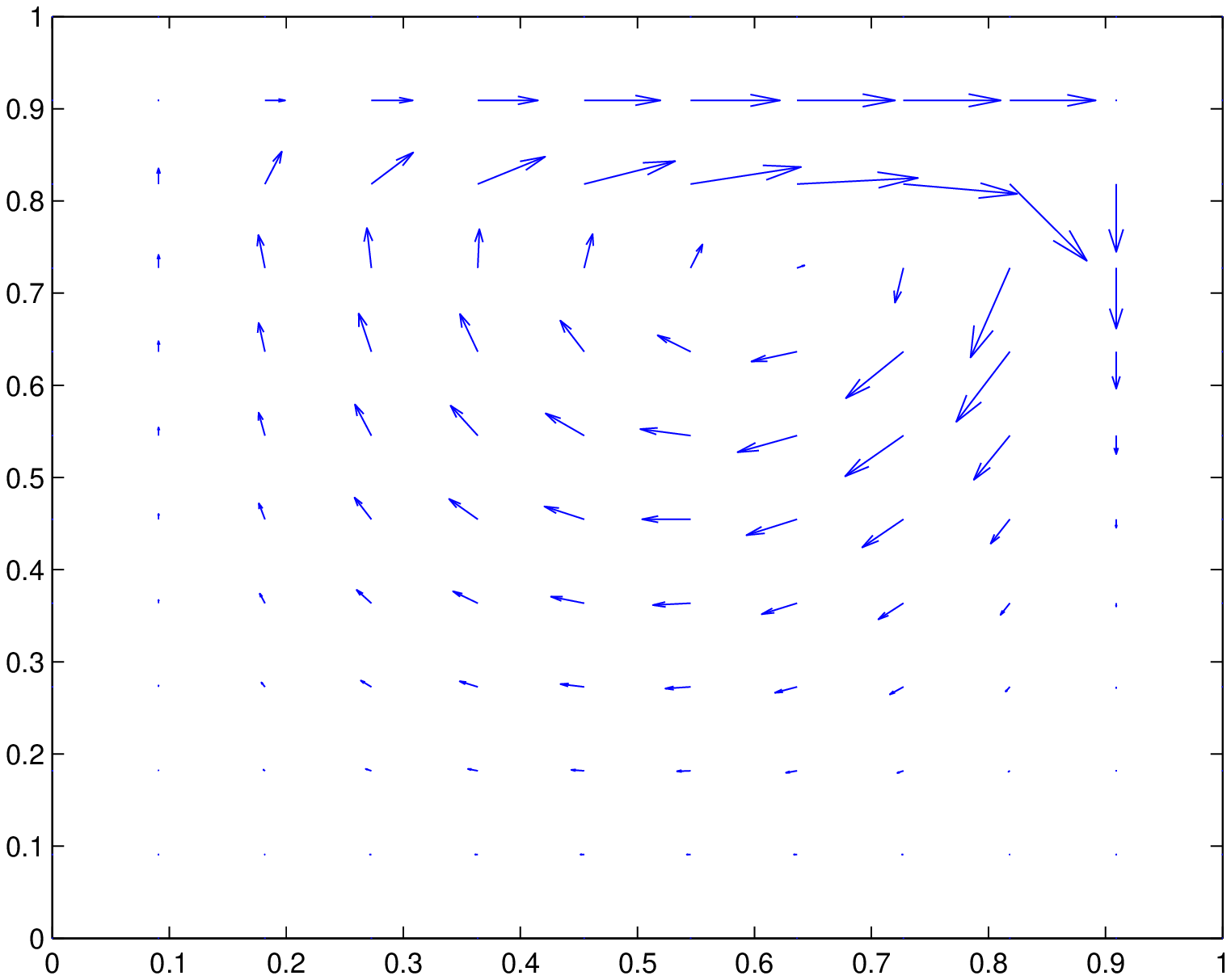}
\includegraphics[width=0.32\textwidth, height=0.16\textheight]{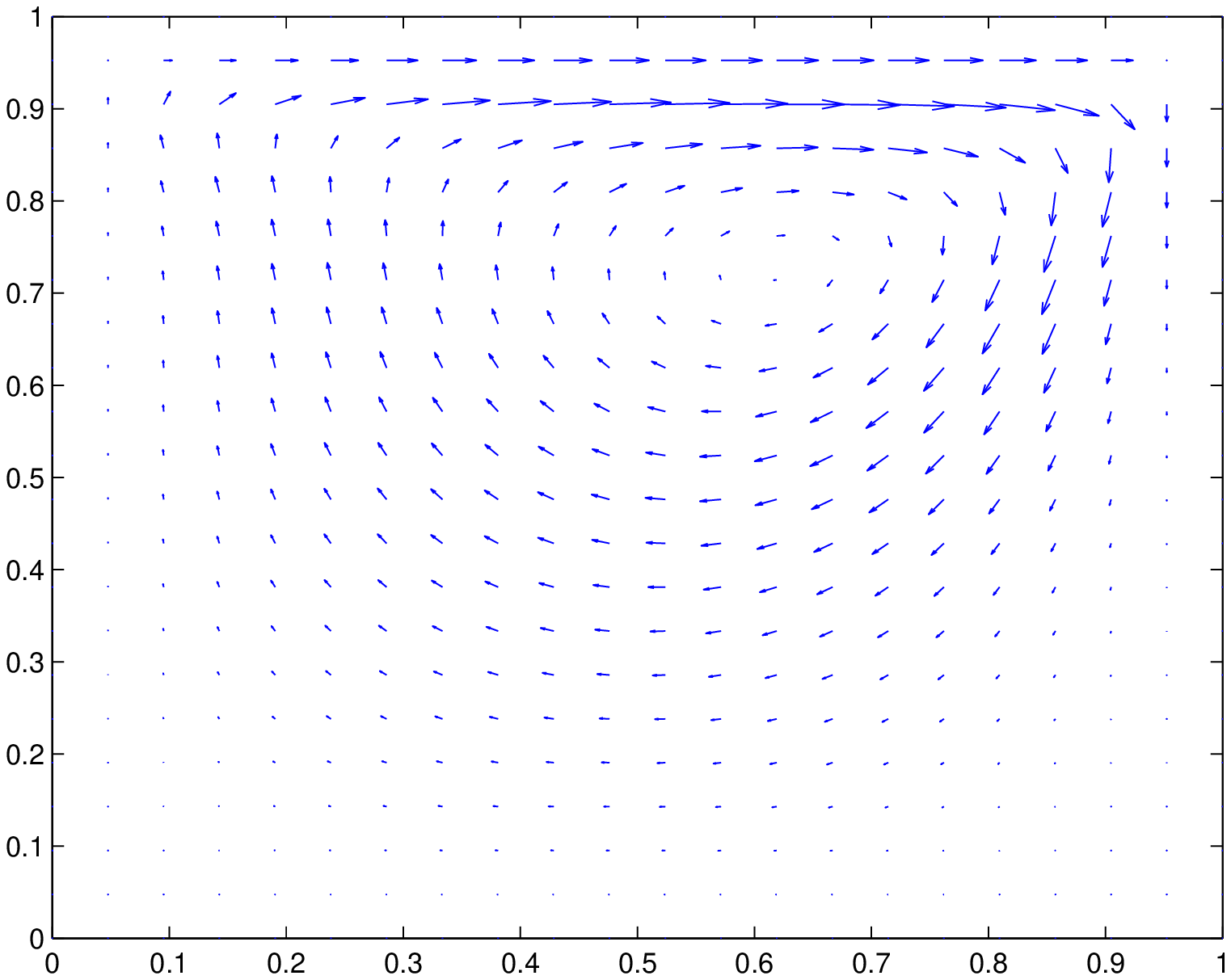}
\includegraphics[width=0.32\textwidth, height=0.16\textheight]{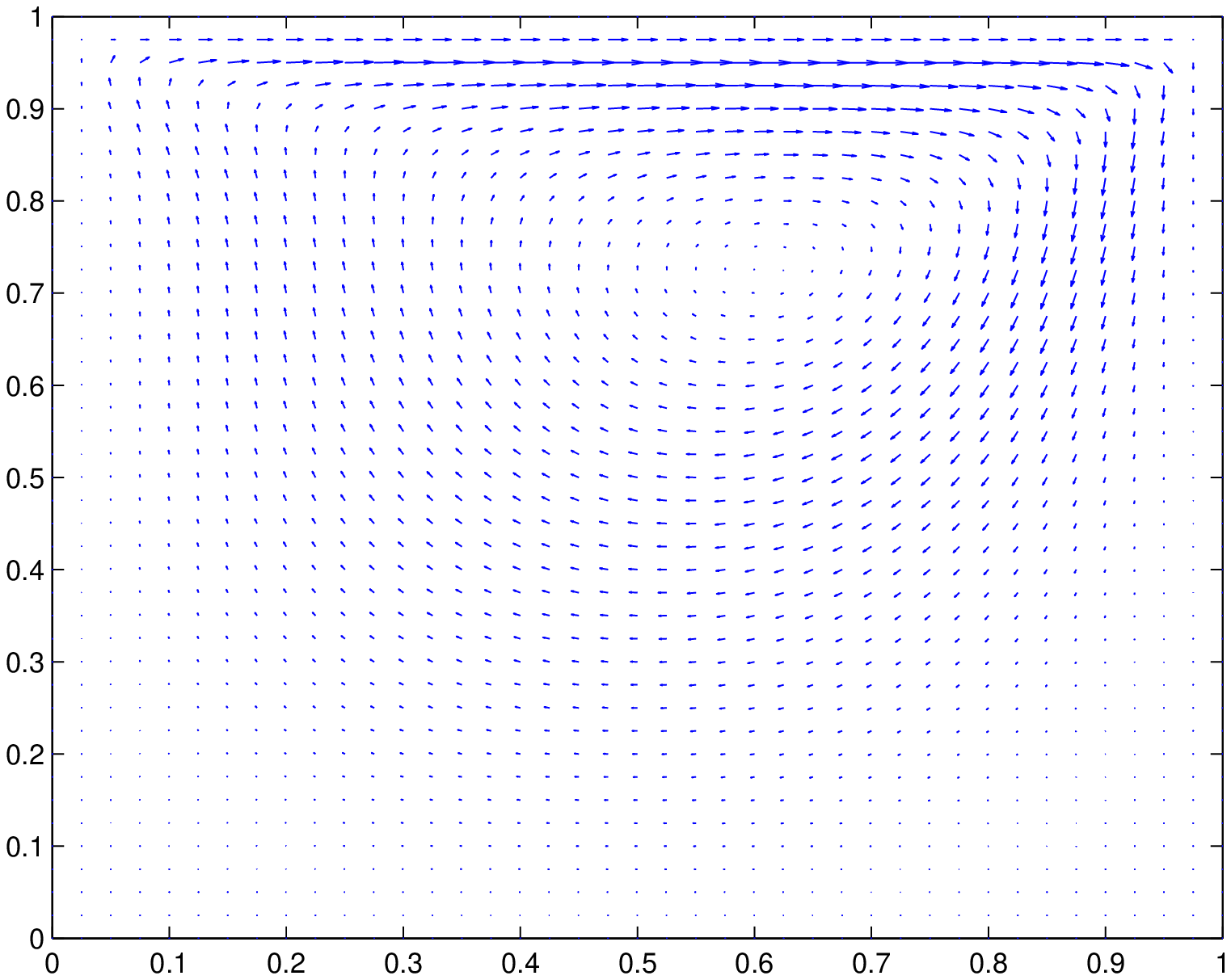}
\caption{$CF(100,1,N)$ solutions for $N=10$ (left), $N=20$ (center), $N=40$ (right)}
\label{Cav100gridreffig}
\end{center}
\end{figure}

We apply SDPR($1$) to $CF(100,1,N)$ and take the solution as starting point for Newton's method. Accurate solutions to the discrete steady cavity flow problem are obtained for $N\in\left\{ 10,\, 15,\, 20 \right\}$. By applying standard grid-refinement methods as in \cite{mknt}, we succeed in extending the solutions to grids of size $30\times 30$ and $40\times 40$. The numerical results are enlisted in Table \ref{Cav100gridrefres} and pictured in Figure \ref{Cav100gridreffig}. Thus, it seems reasonable to conclude, that the minimum energy solution converges to an analytical solution of the steady cavity flow problem. The discrete steady cavity flow problem has multiple solutions. It is an advantage of our method to detect  the minimum kinetic energy solution $u^{\star}$ converging to an analytic solution for $N\rightarrow\infty$.

\begin{table}[ht]
\begin{center}
\begin{tabular}{|r|r|r|r|r|}
\hline $N$ & $w$ & $\epsilon_{\text{sc}}$ &  $t_C$ & $F(u^{\star})$\\
\hline 10& 1 & 4e-11 & 14 & 0.0169\\
\hline 15 & 1 & 6e-16 & 255 & 0.0313\\
\hline 20 & 1 & 6e-16 & 948 & 0.0409\\
\hline 30 & 1 & 4e-11 & 1759 & 0.0503\\
\hline 40 & 1 & 4e-11 & 4156 & 0.0554\\
\hline
\end{tabular}
\caption{Results for $CF(100,1,N)$ for increasing $N$}
\label{Cav100gridrefres}
\end{center}
\end{table}

\subsubsection{Large Reynolds number $R$}
\begin{figure}[ht]
\begin{center}
\includegraphics[width=0.32\textwidth, height=0.16\textheight]{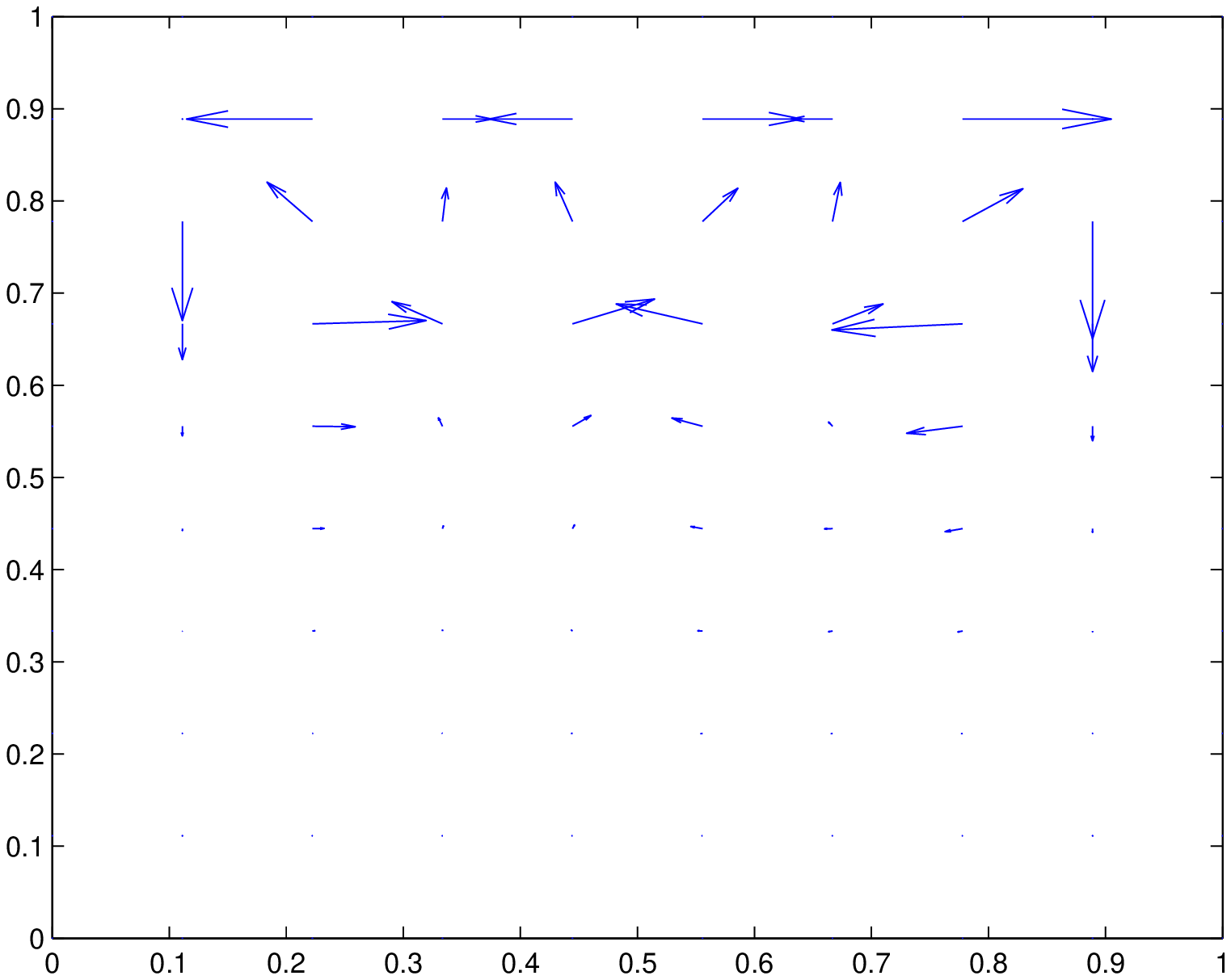}
\includegraphics[width=0.32\textwidth, height=0.16\textheight]{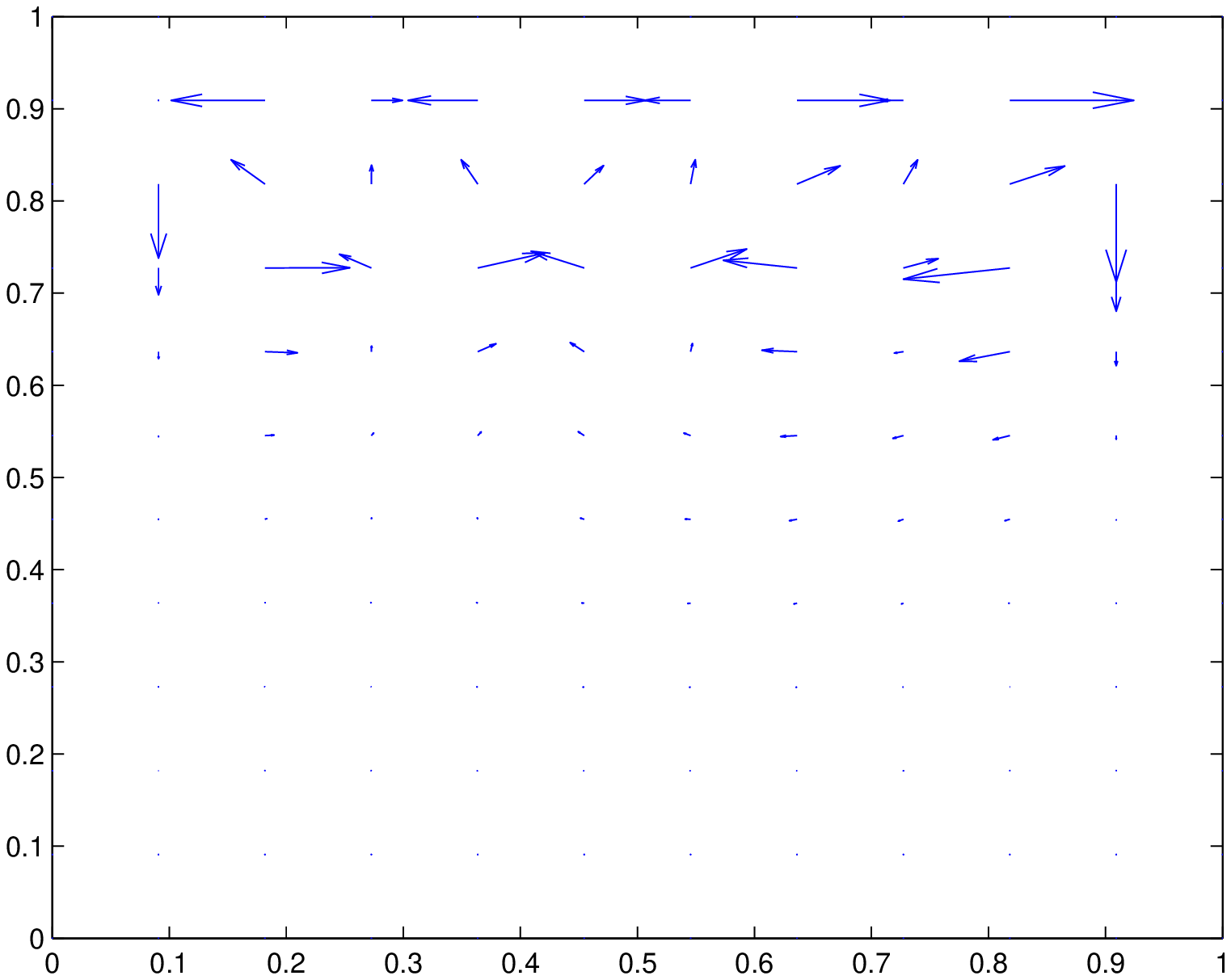}
\includegraphics[width=0.32\textwidth, height=0.16\textheight]{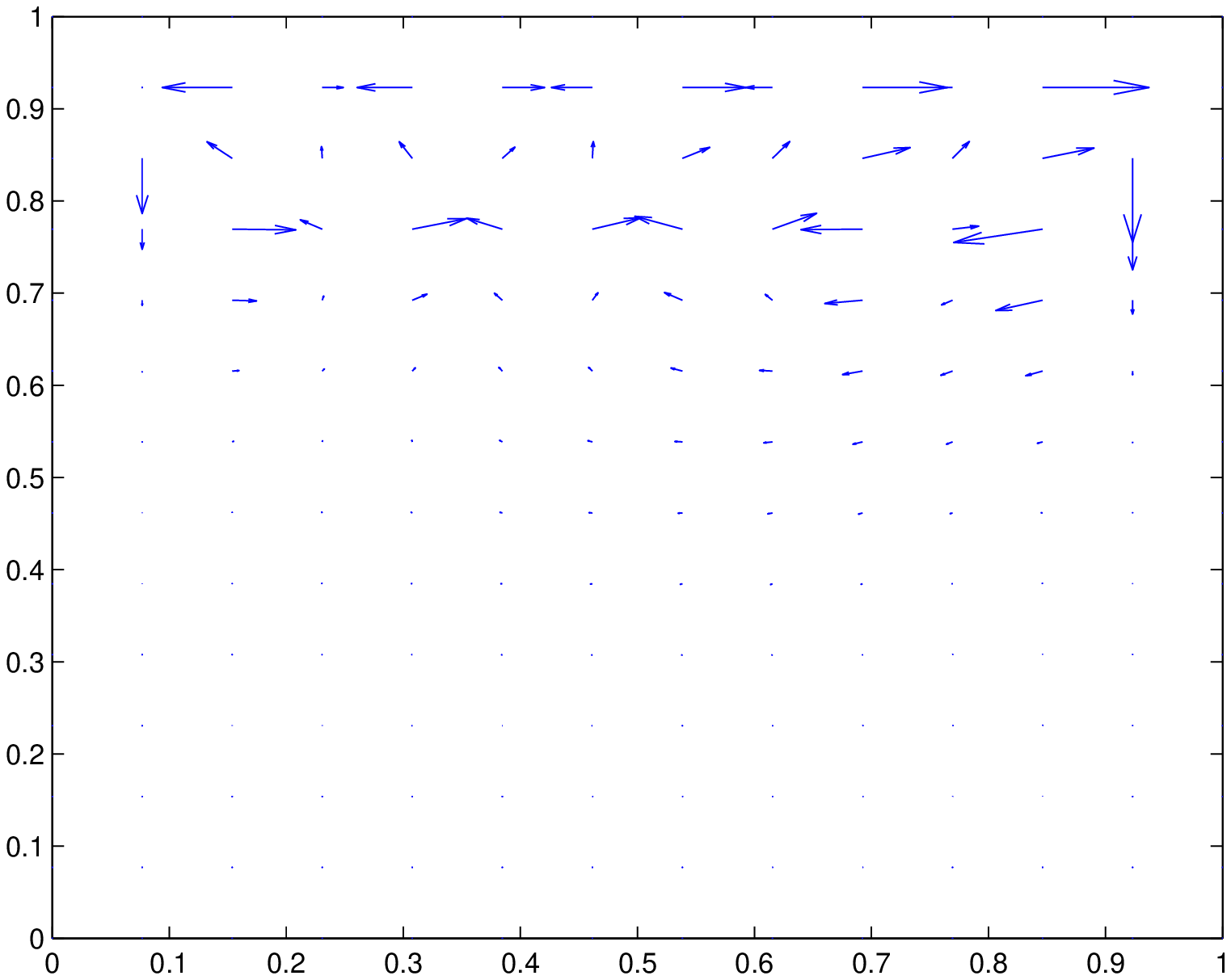}\\
\includegraphics[width=0.32\textwidth, height=0.16\textheight]{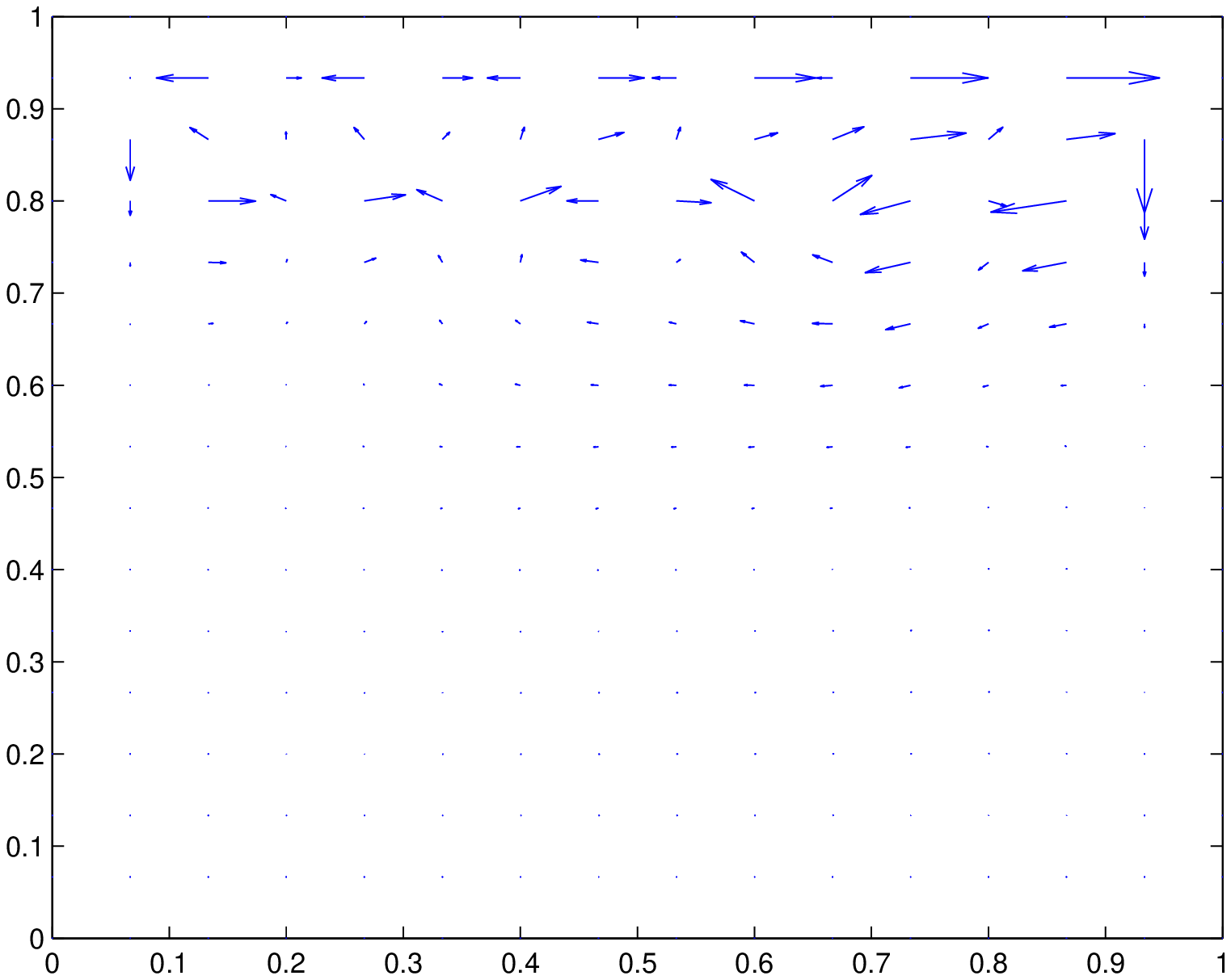}
\includegraphics[width=0.32\textwidth, height=0.16\textheight]{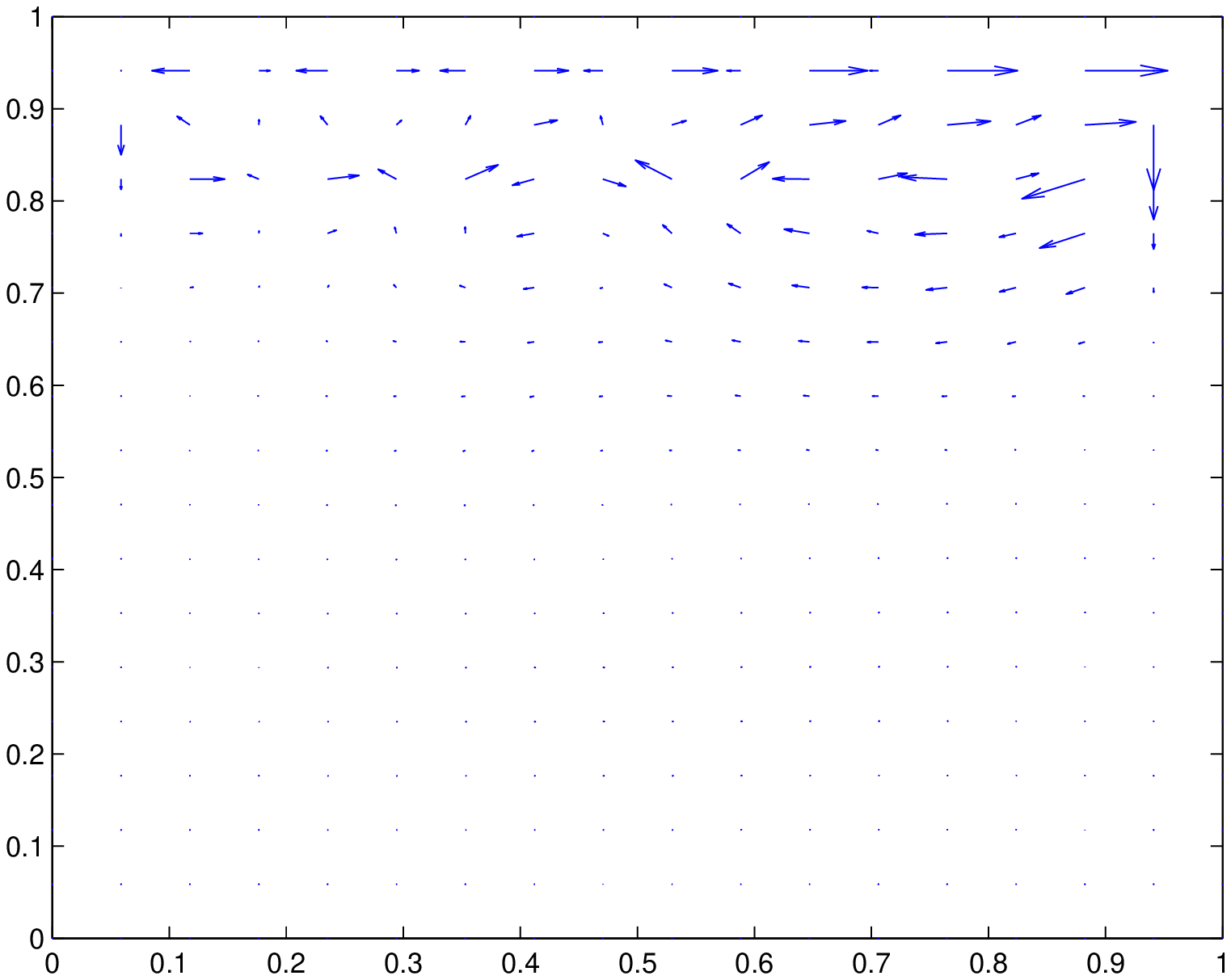}
\includegraphics[width=0.32\textwidth, height=0.16\textheight]{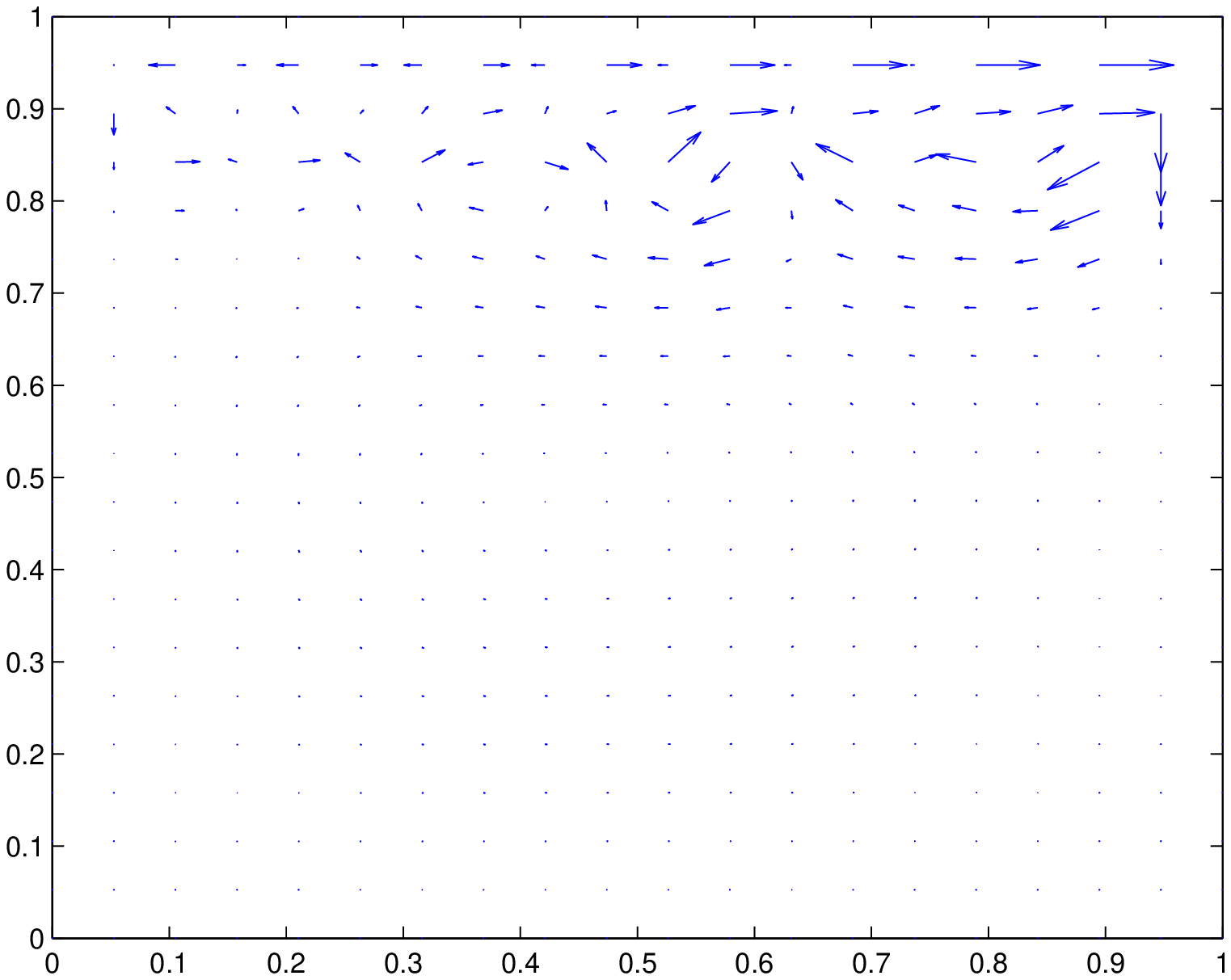}
\caption{SDPR(1) solutions of $CF(10000,1,N)$ for $N=8$ (left, top), $N=10$ (center, top), $N=12$ (right, top), $N=14$ (left, bottom), $N=16$ (center, bottom) and $N=18$ (right, bottom)}
\label{Cav10000Nlargefig}
\end{center}
\end{figure}

In the following we examine $CF(10000,1,N)$ for $N\in\left\{8,\ldots,18\right\}$. For all tested discretizations we were able to find accurate solutions by SDPR($1$) and additionaly applying sequential quadratic programming (SQP). Our results are summarized in Table \ref{Cav10000Nlargeres} and pictured in Figure \ref{Cav10000Nlargefig}.

\begin{table}[ht]
\begin{center}
\begin{tabular}{|r|r|r|r|r|}
\hline $N$ & $w$ & $\epsilon_{\text{sc}}$ &  $t_C$ & $F(u^{(k)})$\\
\hline 8 & 1 & 2e-7 & 7 & 1.5e-6\\
\hline 10& 1 & 3e-10 & 21 & 3.2e-6\\
\hline 12 & 1 & 1e-7 & 49 & 6.0e-6\\
\hline 14 & 1 & 5e-9 & 99 & 1.1e-5\\
\hline 16 & 1 & 4e-12 & 199 & 1.9e-5\\
\hline 18 & 1 & 2e-8 & 501 & 3.9e-5\\
\hline
\end{tabular}
\caption{Results for $CF(10000,1,N)$ for various $N$}
\label{Cav10000Nlargeres}
\end{center}
\end{table}

If we compare the pictures in Figure \ref{Cav10000Nlargefig}, it seems the SDPR(1) solution of $CF(10000,1,N)$ evolves into some {\bf stream-like} solution. Nevertheless, unlike the solutions of $CF(100,1,N)$, we have not been able to expand this solution to a grid of higher resolution by standard interpolation and grid-refinement methods so far. It is possible the solution pictured in Figure \ref{Cav10000Nlargefig} is a fake solution.  

\subsubsection{Large boundary velocity $v$}

\begin{figure}[ht]
\begin{center}
\includegraphics[width=0.32\textwidth, height=0.16\textheight]{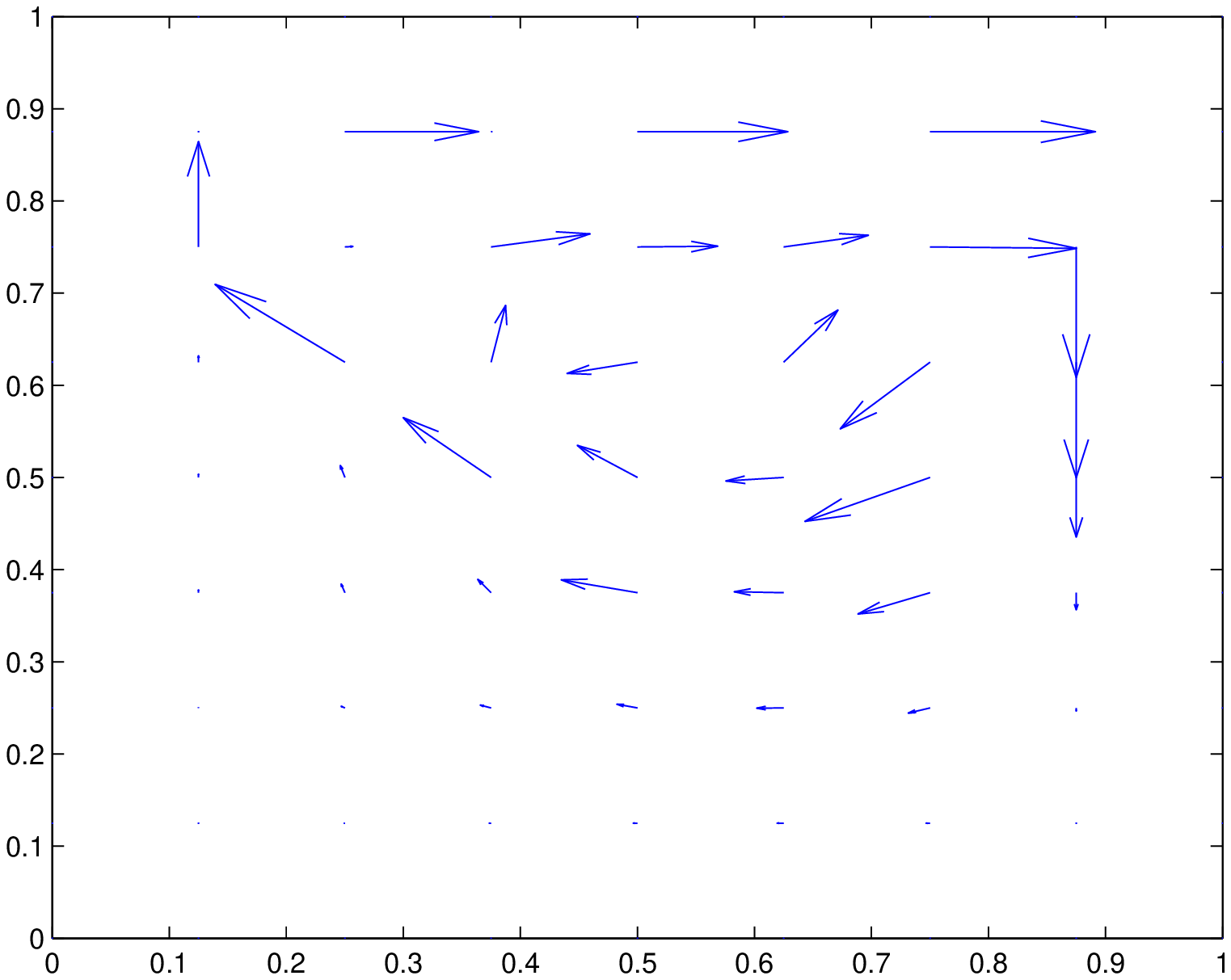}
\includegraphics[width=0.32\textwidth, height=0.16\textheight]{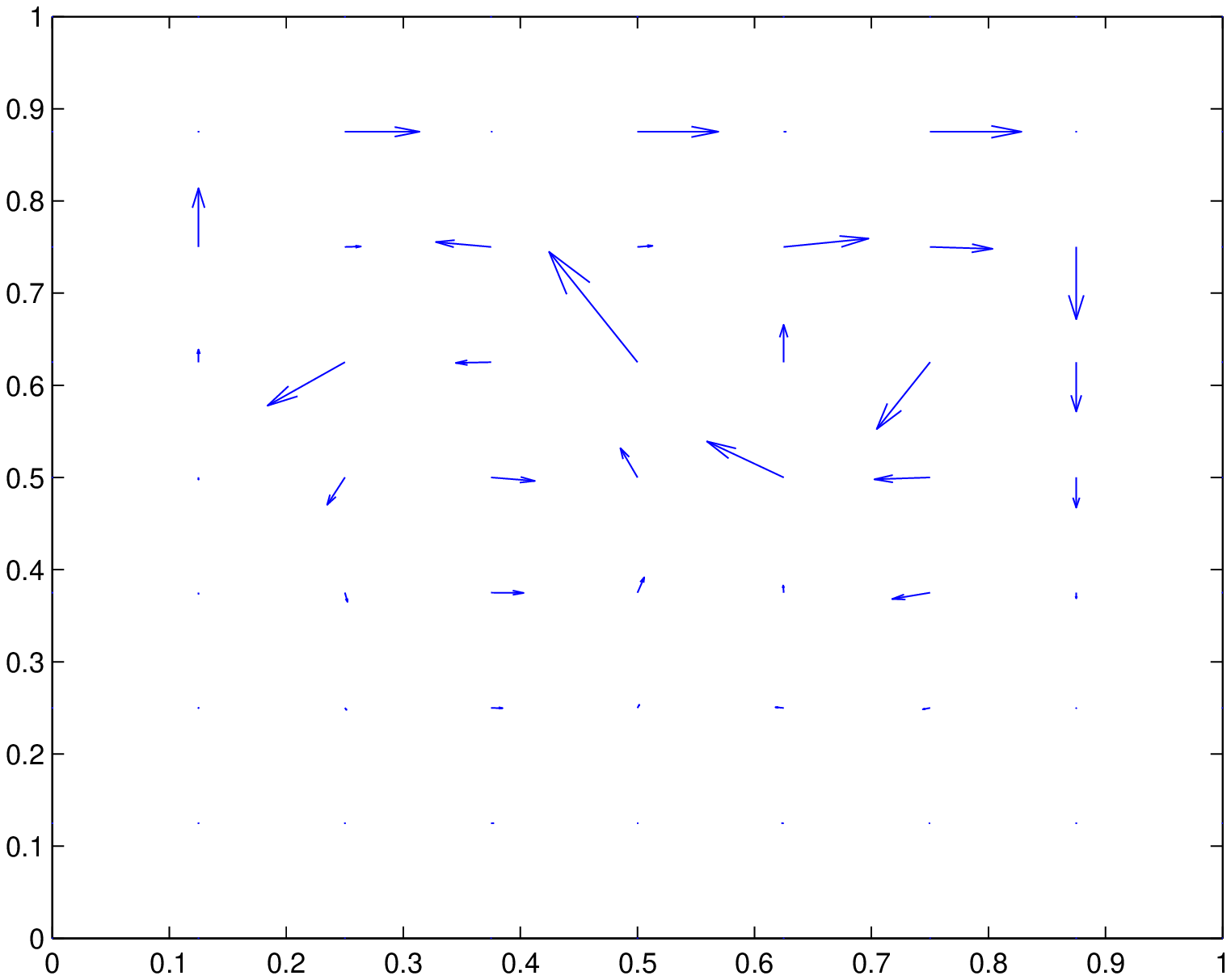}
\caption{SDPR(2) solution (left) and continuation solution (right) for CF(500,10,7)}
\label{Cav500v10fig}
\end{center}
\end{figure}

As an example of a setting with larger boundary velocity we study the problem CF($500,10,7$). We apply the SDPR method with relaxation order $w=2$ and the continuation method, which is a standard method to solve the DSCF we describe in section \ref{secWithContMeth}, and obtain two different solutions for DSCF($500,10,7$), c.f. Table \ref{Cav500v10res} and Figure \ref{Cav500v10fig}. It is interesting to observe that both solutions look like stream solutions: The SDPR(2) solution with one vortex and the continuation solution with two vortices. But as in the previous setting of large $R$, grid-refinement methods to extend these two solutions to higher resolution grids fail. Therefore, it  seems reasonable to conclude that both solutions are fake solutions. Another question is, whether we can derive the continuation solution by Algorithm \ref{enumSolAlg}. As for finding the minimum energy solution, we choose $w=2$. Choosing $b_1^1\geq 2$  in Algorithm \ref{enumSolAlg} generates an SDP that is too large to be solved by the SDP solver SeDuMi. The SDP relaxation for $b_1^1=1$ is tractable, but it is too weak to yield a solution different than the minimum energy solution for various choices of $\epsilon_1^1$.

\begin{table}[ht]
\begin{center}
\begin{tabular}{|l|r|r|r|}
\hline method & $\epsilon_{\text{sc}}$ &  $t_C$ & $F(u^{(k)})$\\
\hline SDPR(2) & 8e-15 & 6774 & 0.0385\\
\hline continuation & 4e-13 & - & 0.0659\\
\hline
\end{tabular}
\caption{Results for $CF(500,10,7)$}
\label{Cav500v10res}
\end{center}
\end{table}

\begin{quest}
Does the minimum kinetic energy solution $u^\star$ of $CF(R,v,N)$ converge to an analytic solution of the steady cavity flow problem for $N\rightarrow\infty$, even for large values of $R$ and $v$?
\label{largeRmotiv}
\end{quest}

\subsubsection{Alternative finite difference scheme}

\begin{figure}[ht]
\begin{center}
\includegraphics[width=0.4\textwidth, height=0.22\textheight]{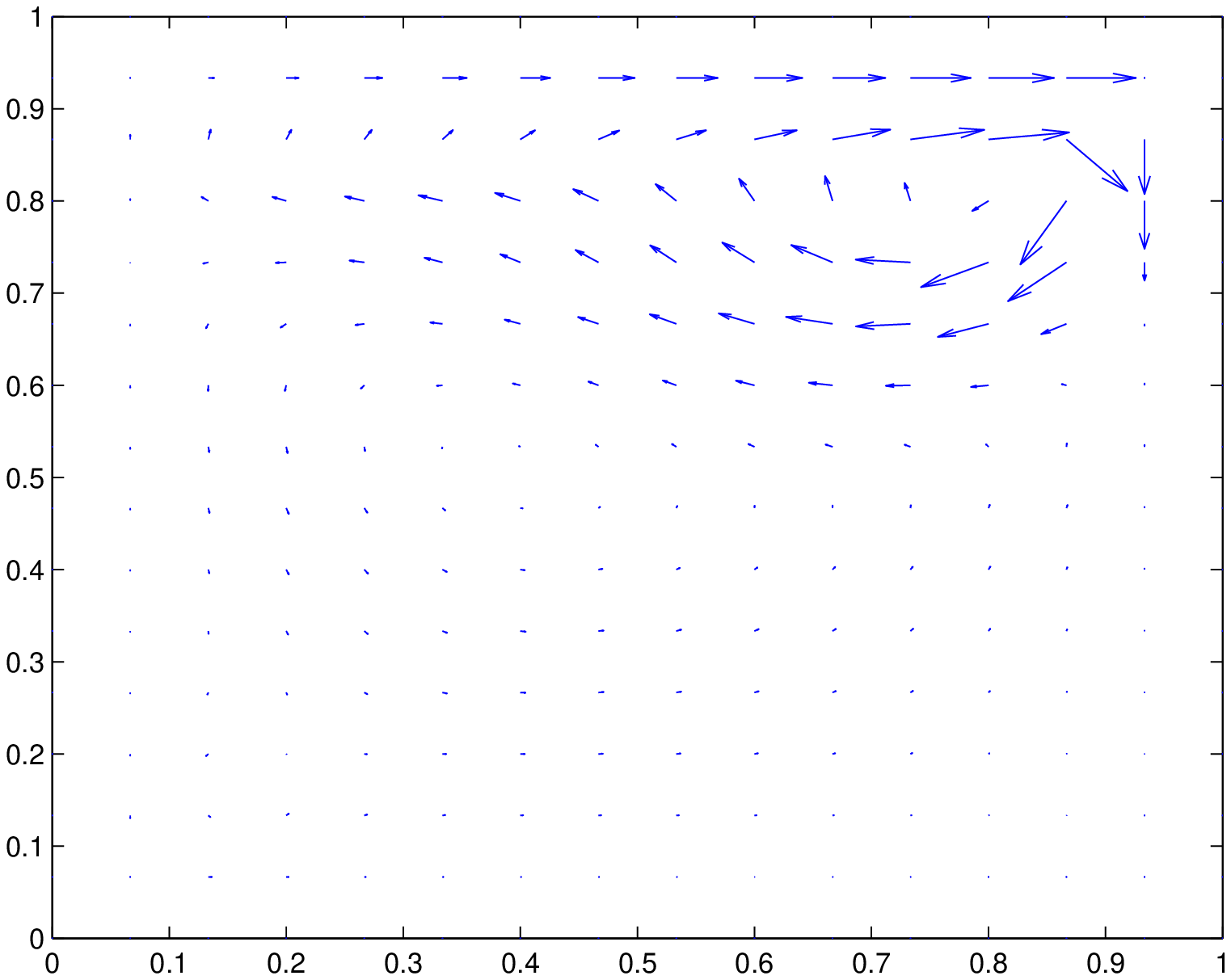}
\includegraphics[width=0.4\textwidth, height=0.22\textheight]{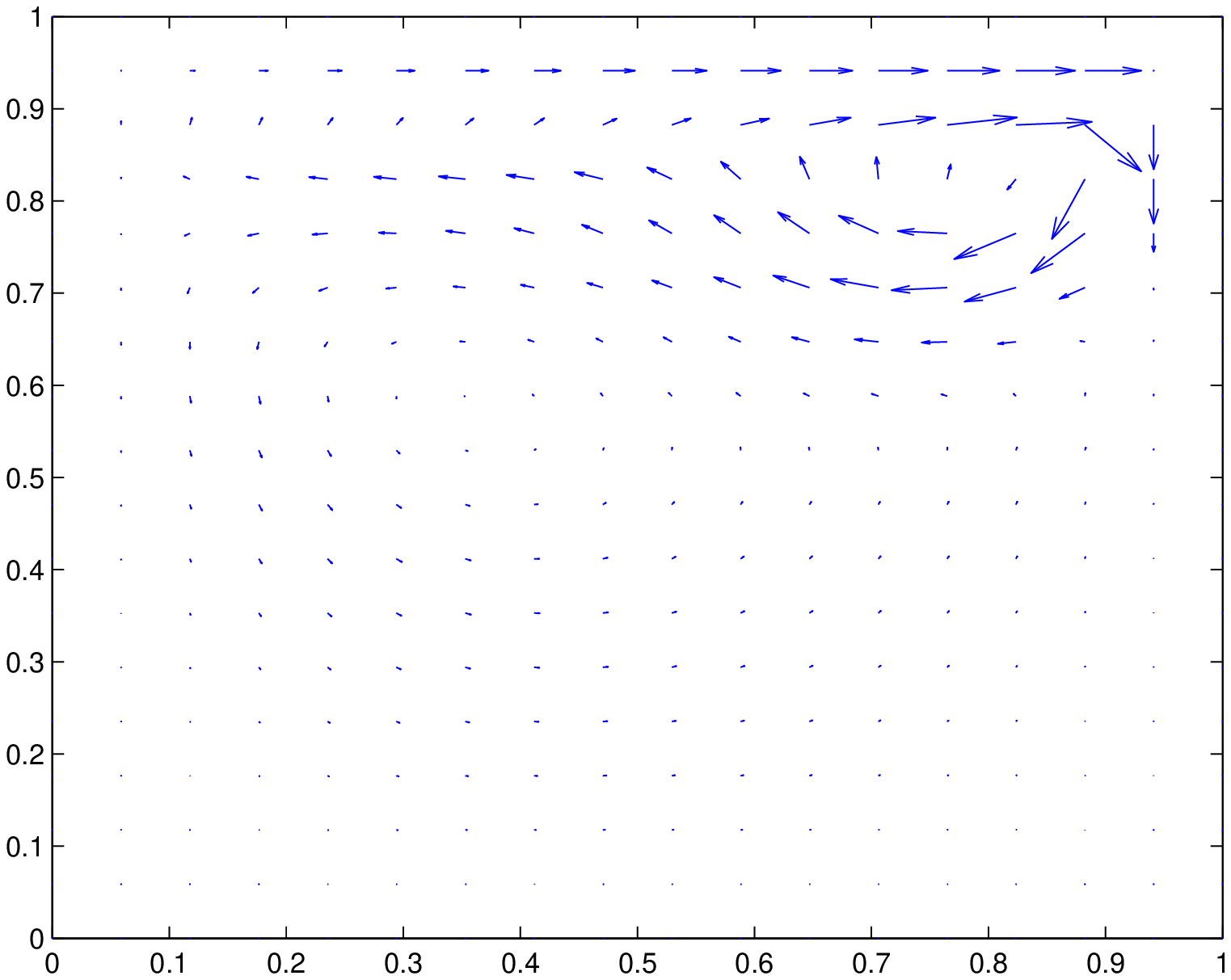}\\
\includegraphics[width=0.4\textwidth, height=0.22\textheight]{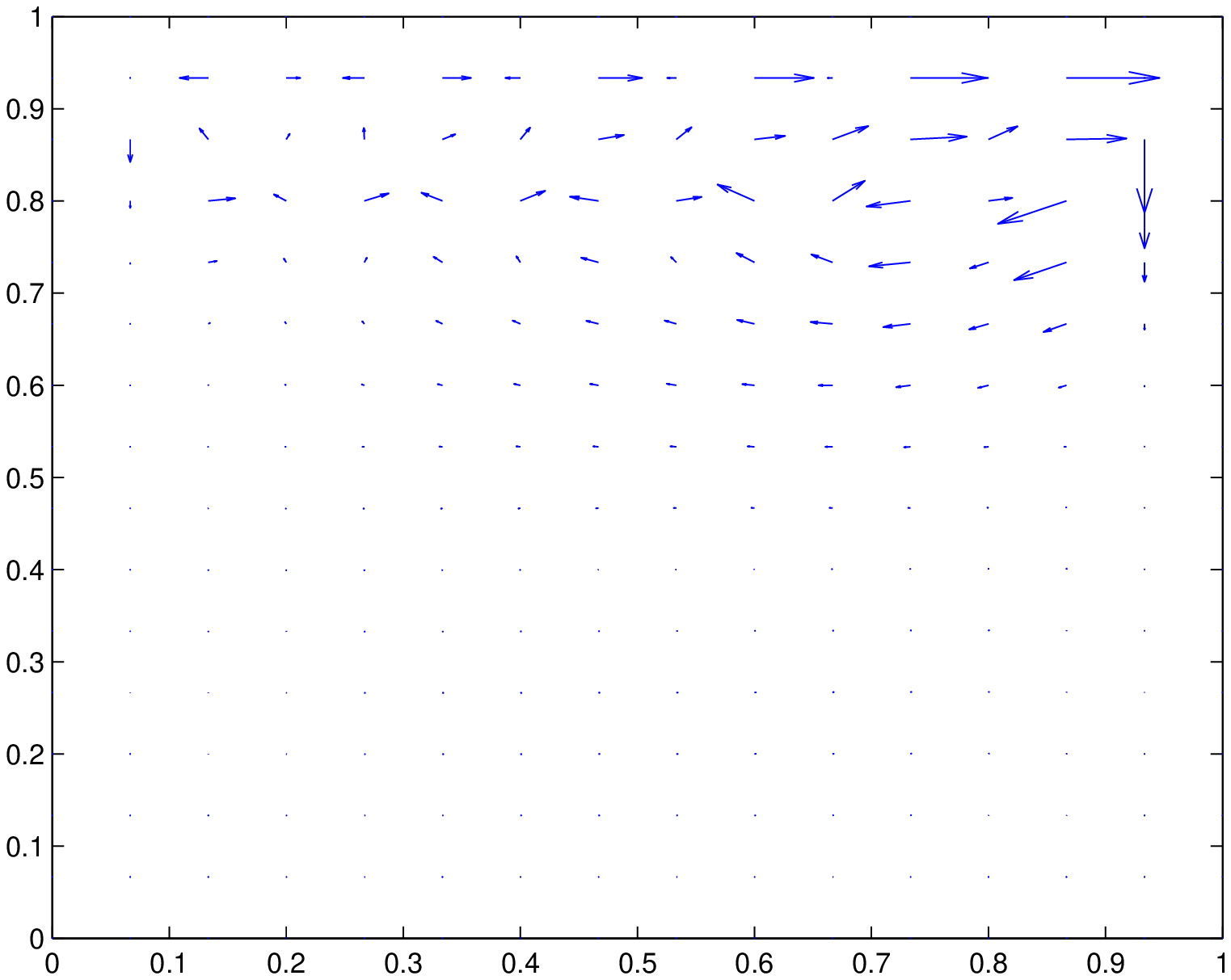}
\includegraphics[width=0.4\textwidth, height=0.22\textheight]{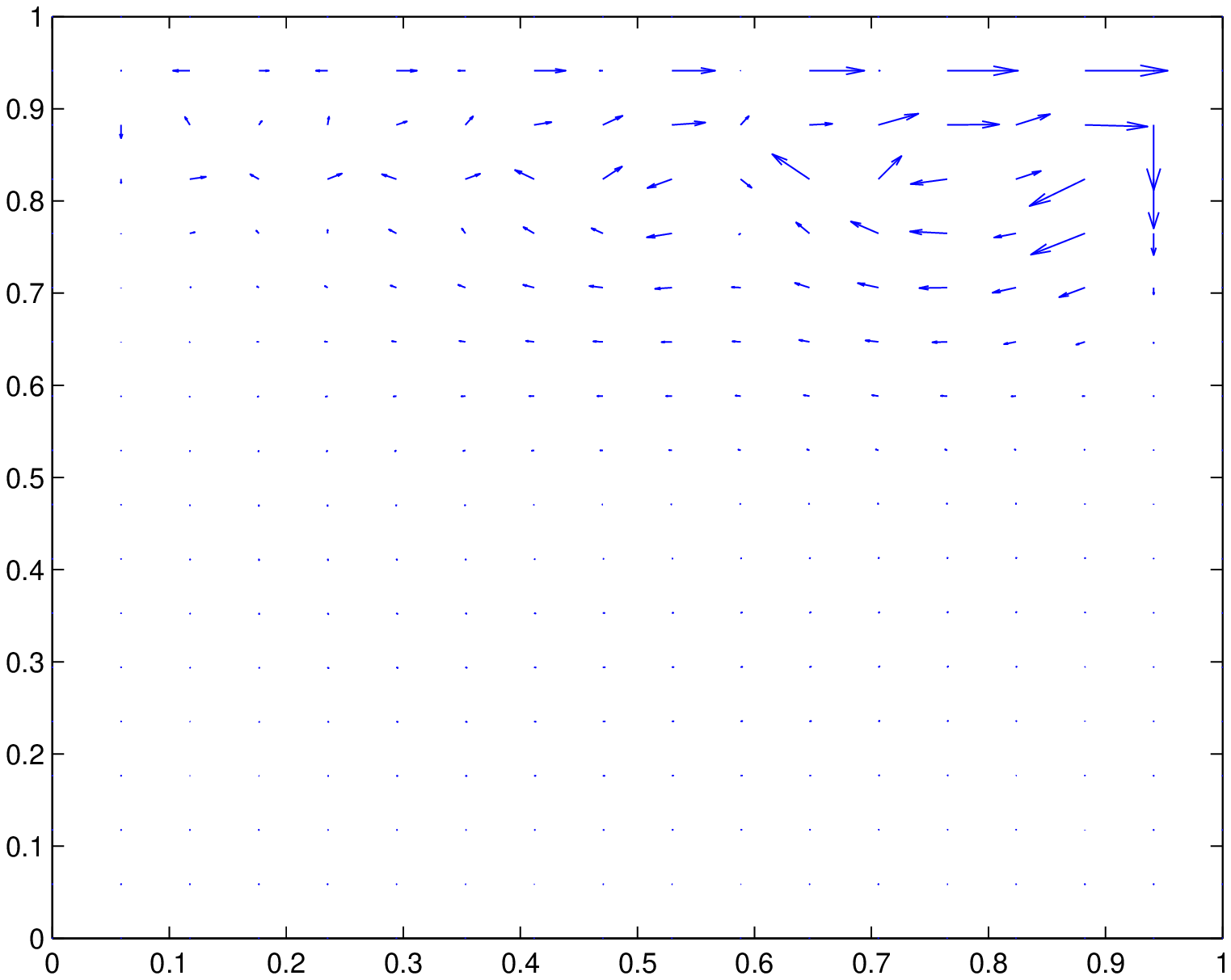}
\caption{Solutions for ADSCF(5000,1,14) (top left), ADSCF(5000,1,16) (top right), DCSF(5000,1,14) (bottom left) and DCSF(5000,1,16) (bottom right)}
\label{figAra}
\end{center}
\end{figure}

We mentioned in Remark \ref{rmkAra}, the simple central finite difference scheme we use does not preserve important physical invariants \cite{arakawa}. Arakawa \cite{arakawa} proposed an alternative finite difference discretization for , that is shown to preserve those invariants. We use this alternative scheme to derive an {\it alternative discrete steady cavity flow problem} ADSCF($R,v,N$) and solve it via the SDPR method. In ADSCF($R,v,N$), the finite difference approximation for $\frac{\partial\psi}{\partial y}\frac{\partial\omega}{\partial x}-\frac{\partial\psi}{\partial x}\frac{\partial\omega}{\partial y}$ in (\ref{eq:omegaconst}) is replaced by
\begin{equation}
\begin{array}{l}
\frac{\partial\psi}{\partial y}\frac{\partial\omega}{\partial x}-\frac{\partial\psi}{\partial x}\frac{\partial\omega}{\partial y}(x_i,\, y_j) \approx\\
-\frac{1}{12h^2} [ \left(\omega_{i,j-1}+\omega_{i+1,j-1}-\omega_{i,j+1}-\omega_{i+1,j+1} \right)\left( \psi_{i+1,j}+\psi_{i,j} \right)\\
 -\left( \omega_{i-1,j-1}+\omega_{i,j-1} - \omega_{i-1,j+1} - \omega_{i,j+1}\right)\left(\psi_{i,j}+\psi_{i-1,j}\right)\\
 +\left(\omega_{i+1,j}+\omega_{i+1,j+1} - \omega_{i-1,j} - \omega_{i-1,j+1}\right) \left( \psi_{i,j+1} + \psi_{i,j} \right)\\
 -\left( \omega_{i+1,j-1} + \omega_{i+1,j} - \omega_{i-1,j-1}-\omega_{i-1,j} \right)\left( \psi_{i,j} + \psi_{i,j-1} \right)\\
 +\left( \omega_{i+1,j}-\omega_{i,j+1} \right)\left( \psi_{i+1,j+1} + \psi_{i,j}\right) - \left(\omega_{i,j-1}-\omega_{i-1,j}\right)\left(\psi_{i,j}+\psi_{i-1,j-1}\right)\\
 +\left( \omega_{i,j+1}-\omega_{i-1,j}\right)\left(\psi_{i-1,j+1}+\psi_{i,j}\right) - \left(\omega_{i+1,j}-\omega_{i,j-1}\right)\left(\psi_{i,j}+\psi_{i+1,j-1}\right)].
\end{array}
\label{araFD}
\end{equation}
It is to be noted that ADSCF($R,v,N$) is less sparse than DSCF($R,v,N$) and it is more difficult to derive accurate solutions by SDPR of relaxation order 1. Nevertheless, we succeed in solving ADSCF($R,v,N$) in some instances. For example, in Table \ref{resAra} and Figure \ref{figAra} we compare the minimum kinetic energy solutions obtained for DSCF($5000,1,N$) and ADSCF($5000,1,N$). It is interesting that the vortex in the minimum kinetic energy solution for ADSCF(5000,1,$N$) is preserved for increasing $N$, whereas the vortex in solution for DSCF(5000,1,$N$) seems to deteriorate.

\begin{table}[ht]
\begin{center}
\begin{tabular}{|r| l |  l |  l |}
\hline Problem & $\epsilon_{\text{sc}}$ & $t_C$ & $F(u^{\star})$\\
\hline ADSCF(5000,1,14) & 7e-12 & 1304 & 1.8e-4\\
\hline ADSCF(5000,1,16) & 5e-10 & 2802 & 3.1e-4\\
\hline DSCF(5000,1,14) & 1e-11 & 419 & 5.6e-4\\
\hline DSCF(5000,1,16) & 3e-10 & 768 & 1.1e-4\\
\hline
\end{tabular}
\caption{Results for solving ADSCF(5000,1,N) compared to DSCF(5000,1,N)}
\label{resAra}
\end{center}
\end{table}

\subsection{Solutions of $CF(R,v,N)$ for increasing Reynolds number $R$ and velocity $v$}
\label{secWithContMeth}

For small Reynolds numbers, we have seen that the minimum kinetic energy solution converges to an analytic solution for $N\rightarrow\infty$ by applying grid-refinement methods. In order to adress Question \ref{largeRmotiv} and to understand why convergence to the analytic solution is a lot more difficult to obtain for large $R$ and $v$, we examine the behavior of the minimum energy solution of the polynomial system DSCF($R,v,N$) and $CF(R,v,N)$, respectively, for increasing Reynolds number $R$.

\subsubsection{Minimum kinetic energy solution for increasing $R$}
To solve the discrete steady cavity flow problem, the proposed SDPR method, Method \ref{sdprm}, is one possibility to find an appropriate starting point for Newton's method. If $w$ is chosen sufficiently large, the output $u$ of Method \ref{sdprm} is guaranteed to accurately approximate the minimum energy solution $u^{\star}$ of $CF(R',v,N)$ and  DSCF($R',v,N$), respectively. In order to show the advantage of the SDPR method we compare our results to solutions of  DSCF($R',v,N$) obtained by a second procedure: In case $R=0$ the discrete steady cavity flow problem DSCF($R,v,N$) is a system of linear equalities, which has an unique solution $u_0(v,N)$, or short $u_0$. Beside solving the linear system, one way to obtain  this solution is solving $CF(0,v,N)$, which is equivalent to solving an SDP as pointed out in Remark \ref{r0rmk}.

\begin{figure}[ht]
\begin{center}
\includegraphics[width=0.49\textwidth, height=0.22\textheight]{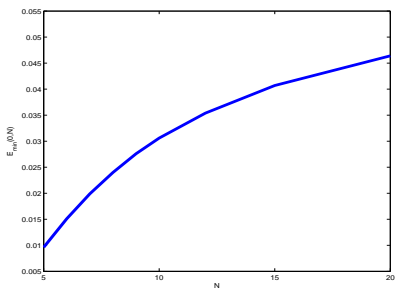} 
\caption{$E_{\min}(0,N)$ for $N\in\left[5, 20\right]$}
\label{Emin0N}
\end{center}
\end{figure}

\begin{method} Naive homotopy-like {\bf continuation method}
\begin{enumerate}
\item
Choose a boundary velocity $v$, grid discretization $N$, Reynolds number $R'$ and step size $\Delta R$.
\item
Solve  DSCF($0,v,N$), i.e. a linear system, and obtain $u^0$. 
\item
Increase $R^{k-1}$ by $\Delta R$: $R^k = R^{k-1} + \Delta R$
\item
Apply Newton's method to DSCF($R^k,v,N$) starting from $u^{k-1}$. Obtain solution $u^k$ as an approximation to the discrete cavity flow problem's solution. 
\item
Iterate 3. and 4. until the desired Reynold's number $R'$ is reached.
\end{enumerate}
\label{homcon}
\end{method}
We call Method \ref{homcon} the {\bf continuation method}. In fact, it is one of the standard methods to find a solution for the steady cavity flow problem. Note, the continuation method does not necessarily yield the minimum kinetic energy solution of DSCF($R,v,N$). In all numerical experiments the boundary velocity $v$ is fixed to $v=1$. Let $u^{\star}(R,N)$ denote the global minimizer of $CF(R,1,N)$, the minimum energy $E_{\min}(R,N)$ is given by $E_{\min}(R,N)=F(u^{\star}(R,N))$. Obviously, it holds $E_{\min}(0,N) = F(u_0(N))$, Figure \ref{Emin0N} shows $E_{\min}(0,N)$ for $N$ ranging from 5 to 20.\\

\begin{figure}[ht]
\begin{center}
\includegraphics[width=0.49\textwidth, height=0.22\textheight]{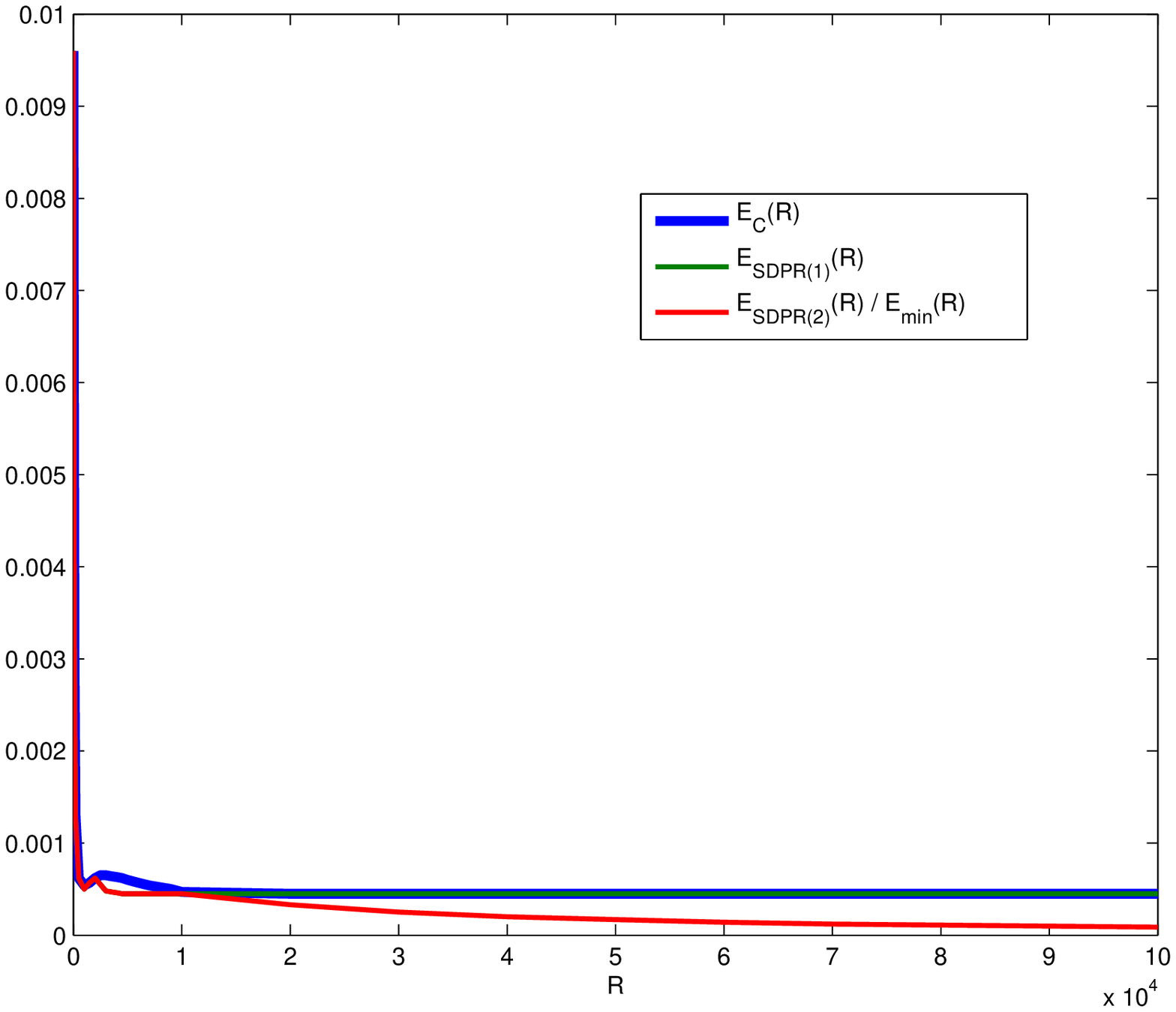}
\includegraphics[width=0.49\textwidth, height=0.22\textheight]{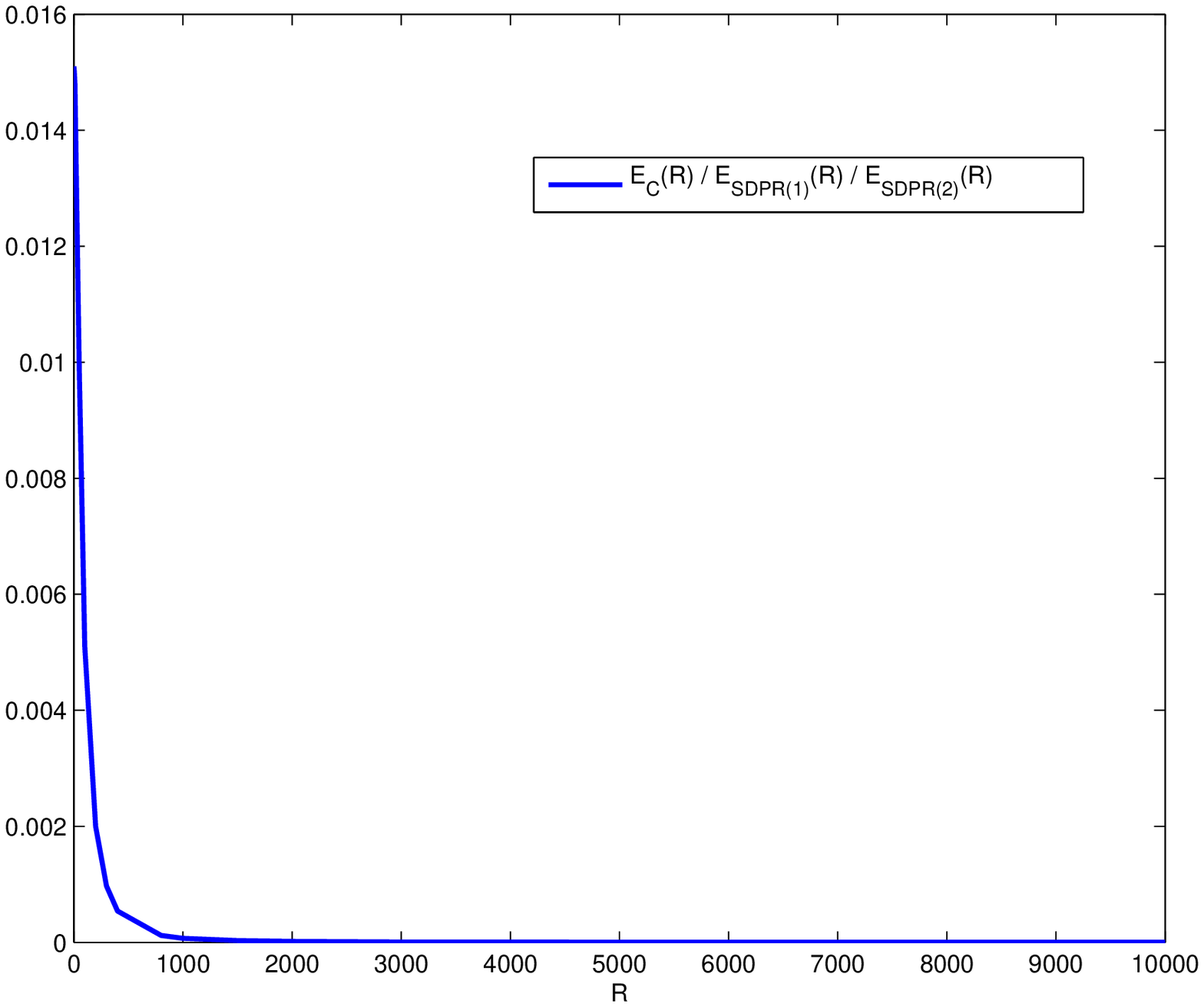}
\caption{$E_C(R)$ and $E_{\min}(R)$ in case $N=5$ (left) and $N=6$ (right)}
\label{ECEminG5R}
\end{center}
\end{figure}

In a next step the solution of DSCF($R,v,N$) obtained by the continuation method starting from $u_0$ is denoted as $\tilde{u}(R)$, and its energy as $E_C(R,N):=F(\tilde{u}(R,N))$.
As illustrated for $N=5$ and $N=7$ in Figure \ref{ECEminG5R}, it is possible for all $R$ to find a continuation $\tilde{u}$ of $u_0$. For $N=5$ the dimension $n$ of the discrete steady cavity flow problem is $n=18$. This dimension is small enough to solve a polynomial system by Gr\"obner basis or polyhedral homotopy method and to determine all complex solutions of the system. Therefore, we can verify whether SDPR combined with SQP detects the global minimizer of $CF(R,v,N)$ or not. It is worth pointing out, that we are able to find the minimum energy solution of the  $CF(R,v,N)$ by applying the SDP relaxation method, whereas this solution cannot be obtained by the standard continuation method. We observe SDPR(1) is sufficient to detect the global optimizer for $R\leq10000$, and for $R\geq 20000$ the global optimizer is obtained by SDPR(2), which is reported in Table \ref{CavRincrG5}. 

\begin{table}[ht]
\begin{center}
\begin{tabular}{|r|r|r|r|r|r|r|}
\hline $R$ & $N_{\CC}$ & $N_{\RR}$ & $E_C(R)$ & $E_{SDPR(1)}$ & $E_{SDPR(2)}$ & $E_{\min}$\\
\hline 0 & 1 & 1 & 0.0096 & 0.0096 & 0.0096 & 0.0096\\
\hline 1 & 35 & 23 & 0.0096 & 0.0096 & 0.0096 & 0.0096\\
\hline 10 & 37 & 17 & 0.0094 & 0.0094 & 0.0094 & 0.0094\\
\hline 100 & 37 & 13 & 0.0030 & 0.0030 & 0.0030 & 0.0030 \\
\hline 200 & 37 & 11 & 0.0013 & 0.0013 & 0.0013 & 0.0013 \\
\hline 500 & 37 & 13 & 6.2e-4 & 6.2e-4 & 6.2e-4 & 6.2e-4\\
\hline 1000 & 37 & 13 & 5.4e-4 & 5e-4 & 5e-4 & 5e-4\\
\hline 2000 & 37 & 13 & 6.2e-4 & 6.2e-4 & 6.2e-4 & 6.2e-4\\
\hline 3000 & 38 & 18 & 6.5e-4 & 4.8e-4 & 4.8e-4 & 4.8e-4\\
\hline 4000 & 37 & 17 & 6.3e-4 & 4.6e-4 & 4.6e-4 & 4.6e-4\\
\hline 6000 & 36 & 16 & 5.7e-4 & 4.5e-4 & 4.5e-4 & 4.5e-4 \\
\hline 8000 & 36 & 16 & 5.2e-4 & 4.5e-4 & 4.5e-4 & 4.5e-4 \\
\hline 10000 & 35 & 17 & 4.7e-4 & 4.5e-4 & 4.5e-4 & 4.5e-4\\
\hline 20000 & 35 & 17 & 4.5e-4 & 4.5e-4 & 3.3e-4  & 3.3e-4\\
\hline 30000 & 35 & 17 & 4.5e-4 & 4.5e-4 & 2.5e-4 & 2.5e-4\\
\hline 50000 & 35 & 17 & 4.5e-4 & 4.5e-4 & 1.7e-4 & 1.7e-4\\
\hline 70000 & 35 & 16 & 4.5e-4 & 4.5e-4 & 1.2e-4 & 1.2e-4\\
\hline 100000 & 34 & 16 & 4.5e-4 & 4.5e-4 & 8.8e-5 & 8.8e-5\\
\hline 
\end{tabular}
\caption{Numerical results for $CF(R,1,5)$}
\label{CavRincrG5}
\end{center}
\end{table}

\begin{figure}[ht]
\begin{center}
\includegraphics[width=0.49\textwidth, height=0.22\textheight]{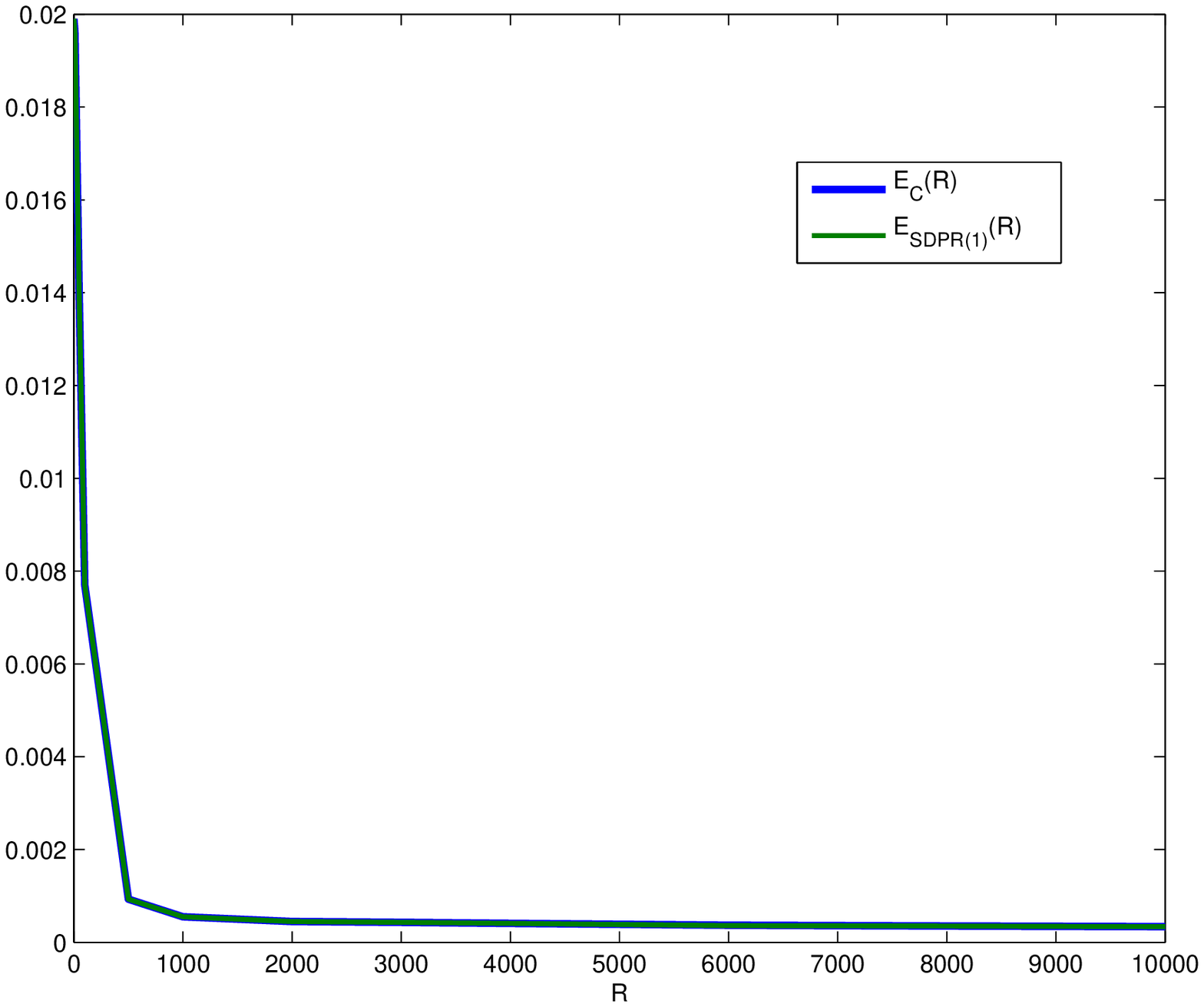}
\caption{$E_C(R)$ and $E_{SDPR(1)}(R)$ in case $N=7$}
\label{ECESG7R}
\end{center}
\end{figure}

In case of $N=6$ and $N=7$ the dimension of the polynomial system is too large to be solved by Gr\"obner basis or polyhedral homotopy method for $R>0$. For $N=6$ the continuation method, SDPR(1) and SDPR(2) yield the same solution for all tested $R$ as pictured in Figure \ref{ECEminG5R} (right). And in case of $N=7$ the continuation solution $\tilde{u}(R)$ is detected by SDPR(1) as well, except the case $R=6000$, where a solution with slightly smaller energy is detected by SDPR(1), as documented in Table \ref{CavRincrG7} and illustrated in Figure \ref{ECESG7R}.

\begin{table}[ht]
\begin{center}
\begin{tabular}{|l|r|r|r|r|r|r|r|r|r|}
\hline $R$ & 0 & 50 & 100 & 500 & 2000 & 4000 & 6000 & 8000 & 10000\\
\hline $E_C(R)$ & 2.0e-2 & 1.4e-2 & 7.7e-3 & 9.3e-4 & 4.5e-4 & 4.1e-4 & 3.7e-4 & 3.5e-4 & 3.4e-4\\
\hline $E_{\text{SDPR}(1)}$ & 2.0e-2 & 1.4e-2 & 7.7e-3 & 9.3e-4 & 4.5e-4 & 4.1e-4 & 3.6e-4 & 3.5e-4 & 3.4e-4\\
\hline  
\end{tabular}
\caption{Numerical results for $CF(R,1,7)$}
\label{CavRincrG7}
\end{center}
\end{table}

Summarizing these results, $F(u_0(R,N))\geq F(\tilde{u}(R,N))$ for any of the tested $R>0$. It is an advantage of our approach to show, $\tilde{u}(R,N)$ is in general not the optimizer of $CF(R,1,N)$ for increasing $R$. In fact, for some settings we obtain far better approximations to the minimum energy solution than $\tilde{u}(R,N)$. Furthermore,  $E_{min}(R)$ and $E_C(R)$ are both decreasing in $R$. The behavior of $E_C,\, E_{SDPR}$ and $E_{\min}$ coincides for all chosen discretization $N$ and motivates the following conjecture.

\begin{conj} Let boundary velocity $v$ and discretization $N$ be fixed.
\begin{enumerate}
\item[a) ]
$F(u_0(v,N)) = E_{\min}(0,v,N) \geq E_{\min}(R,v,N) \geq 0\quad\forall R\geq0$.
\item[b) ]
$E_{\min}(R,v,N) \rightarrow 0$ for $R\rightarrow\infty$.
\end{enumerate}
\label{Rconj}
\end{conj}

As an application, Conjecture \ref{Rconj} can be used as a certificate for the non-optimality of a feasible solution $u'$ of $CF(R,v,N)$ in case $F(u'(R,v,N))> E_{\min}(0,v,N)$. In fact, as it seems to be always possible to extend $u_0$ via continuation method, $\tilde{u}(R,v,N)$ can serve as a non-optimality certificate in case $F(u'(R,v,N)) > F(\tilde{u}(R,v,N))$.

\subsubsection{Stability analysis}
\begin{figure}[ht]
\begin{center}
\includegraphics[width=0.32\textwidth, height=0.16\textheight]{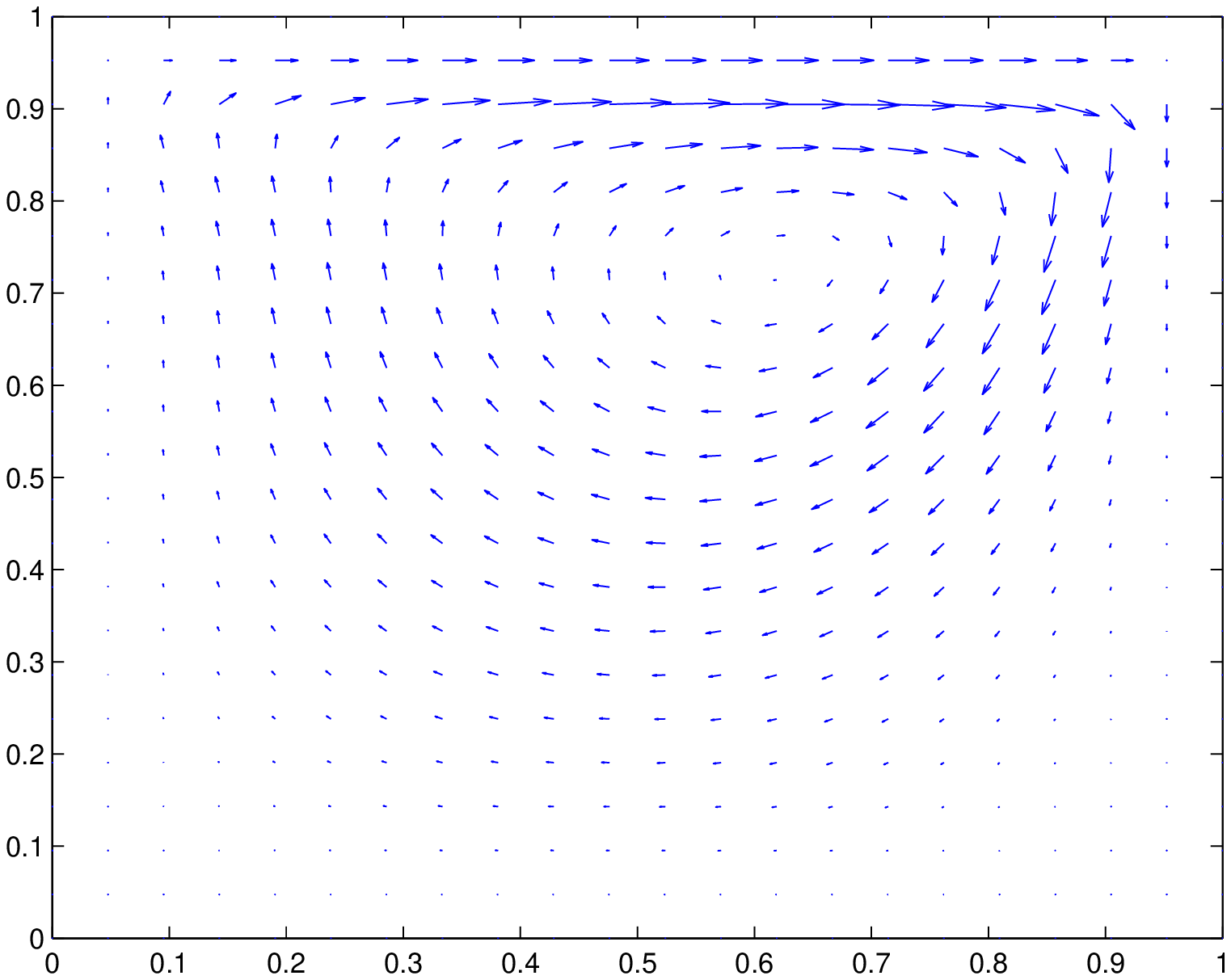}
\includegraphics[width=0.32\textwidth, height=0.16\textheight]{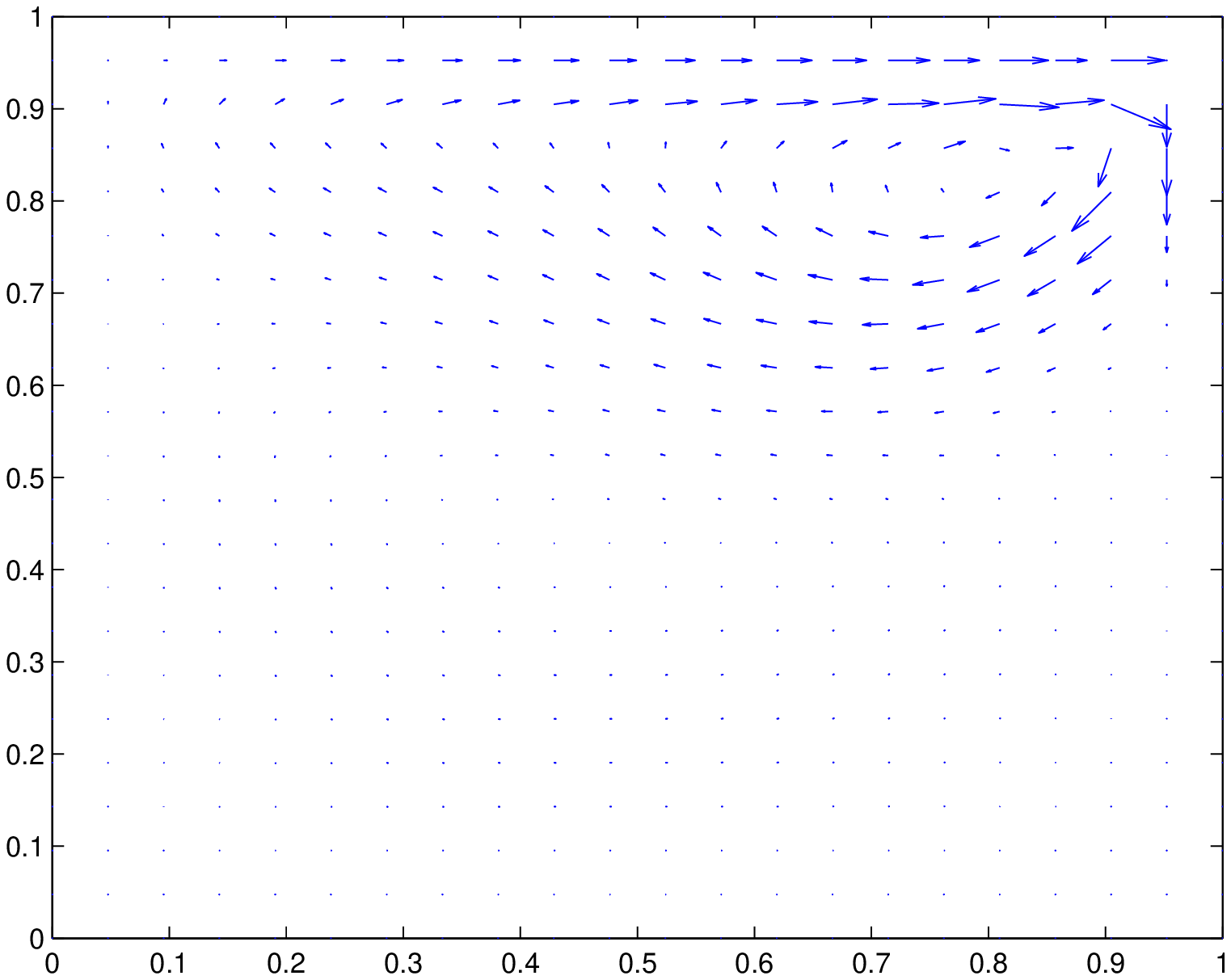}
\includegraphics[width=0.32\textwidth, height=0.16\textheight]{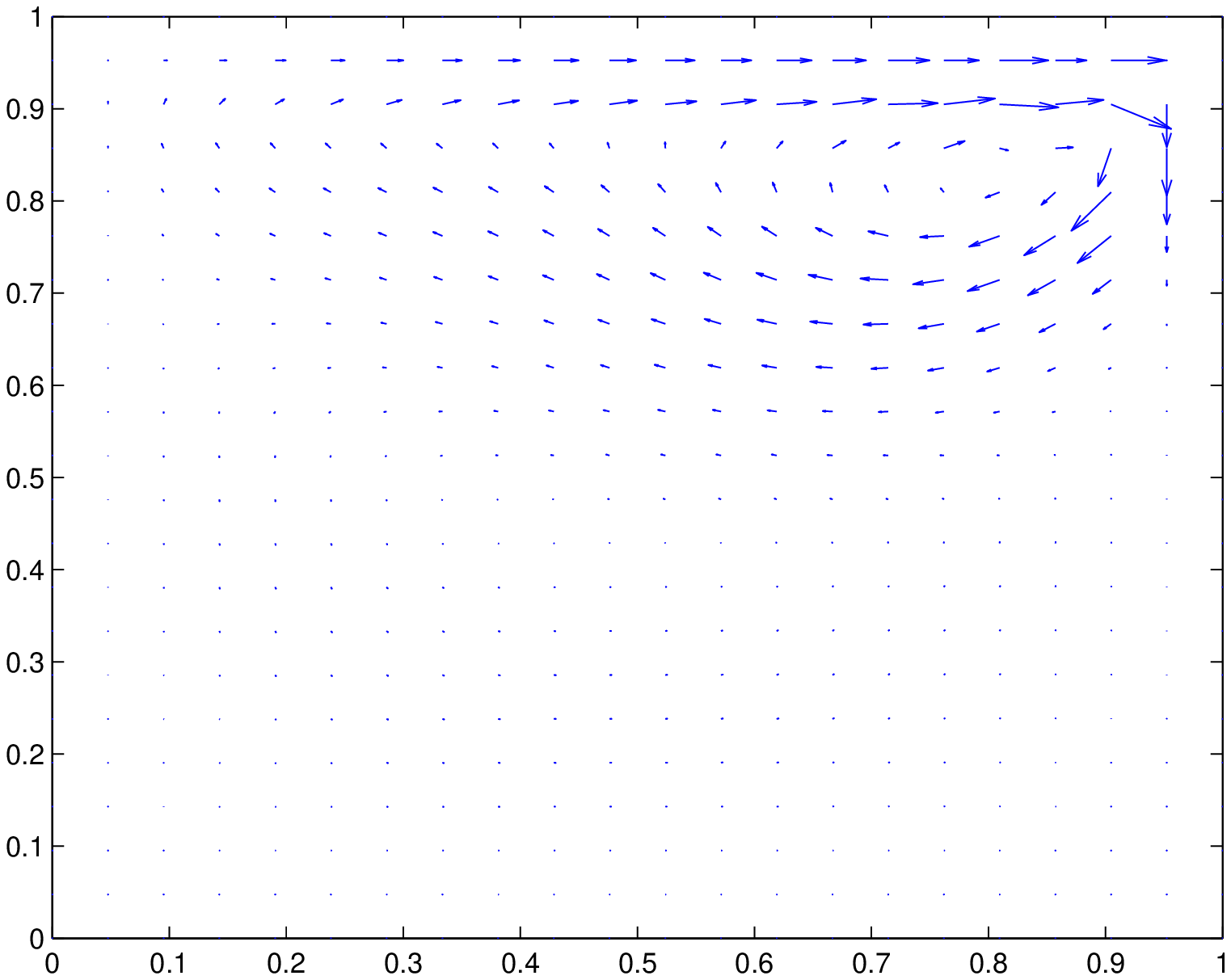}
\caption{$u^{\star}(R)$ for $R=100$ (left), $R=776$ (center) and $R=800$ (right)}
\label{StabanG20fig}
\end{center}
\end{figure}

Finally, we examine the stability of the  minimum kinetic energy solution. We fix $N$ and $v=1$ and increase the Reynolds number $R$. Let $J(u)$ denote the Jacobi matrix of the polynomial system DSCF($R,v,N$) at the solution $u$, $\lambda_{\max}(u)$ its maximal eigenvalue and $N_\lambda^{+}$ its number of positive eigenvalues. A solution $u$ to DSCF($R,v,N$) is called {\it stable}, if all eigenvalues of $J(u)$ are non-negative, otherwise it is called {\it unstable}. In case $N=10$ and $N=20$, we observed that the minimum kinetic energy solution $u^{\star}(R)$ is stable for small $R$, and as  $R$ is increased and exceeds some threshold $R'(v,N)$, $u^{\star}(R)$ becomes unstable. See Table \ref{StabanG20res} and Figure \ref{StabanG20fig}.

\begin{table}[ht]
\begin{center}
\begin{tabular}{|r|r|r|r|r|r|}
\hline $R$ & $N$ & $\epsilon_{\text{scaled}}$ & $t_C$ & $\lambda_{\max}(u)$ & $N_\lambda^{+}$\\
\hline 100 & 20 & 9e-16 & 491 & -0.0636 & 0\\
\hline 750 & 20 & 5e-14 & 366 & -0.0615 & 0 \\
\hline 775 & 20 & 6e-16 & 486 & -0.0014 & 0 \\
\hline 776 & 20 & 4e-16 & 486 & 0.0010 & 2 \\
\hline 800 & 20 & 7e-16 & 486 & 0.0580 & 2\\
\hline 1000 & 20 & 1e-12 & 527 & 0.5124 & 2\\
\hline 300 & 10 & 6e-16 & 10 & -0.3069 & 0\\
\hline 350 & 10 & 6e-16 & 11 & -0.1360 & 0\\
\hline 400 & 10 & 3e-16 & 9 & 0.1217 & 2\\
\hline 500 & 10 & 4e-16 & 9 & 0.6043 & 4\\
\hline
\end{tabular}
\caption{The stable solution $u^{\star}(R)$ becomes unstable for increasing $R$.}
\label{StabanG20res}
\end{center}
\end{table}

\section{Conclusion}

We proposed an algorithm to enumerate all solutions of the discrete cavity flow problem with respect to their kinetic energy. Our algorithm takes advantage of the sparse semidefinite relaxation method (SDPR) in order to find a good starting point for Newton's method or sequential quadratic programming. We can guarantee the convergence of the algorithm's output to the smallest kinetic energy solutions of the polynomial system, if the order of the SDPR tends to infinity.  Our numerical experiments for various choices of $R$ and $v$ have demonstrated that it is sufficient to apply SDPR of order one or two, to succeed in obtaining accurate approximations to the smallest energy solutions of the discrete cavity flow problem by our enumeration algorithm. In case of small Reynolds numbers our algorithm allowed another interesting observation: Among all solutions of the polynomial system given by the discrete cavity flow problem, the minimal kinetic energy solution converges to an analytic solution of the continuous steady cavity flow problem. 
In case of large Reynolds number $R$ we are not able to extend our coarse grid solutions to a finer grid, yet, although many of them look like stream solutions when the kinetic energy is small. 
It is known that the set of solutions of the discrete cavity flow problem contains lots of non-physical solutions or fake solutions, but there has been no systematic study of the discrete cavity flow problem as a polynomial system so far. Moreover, the more interesting stream-like solutions of the discrete steady cavity flow problem are usually among the 3rd or 4th smallest kinetic energy solutions. Our enumeration algorithm based on the SDPR method provides a powerful tool to detect the smallest energy solutions one by one, which is a strong advantage compared to the existing methods. The further analysis of the polynomial system derived from the steady cavity flow problem for large Reynolds number $R$ will remain an interesting topic in future. 
Also, the conjecture that the mininum kinetic energy converges to zero for increasing $R$ is left for future research and may constitute an interesting property of the minimum energy solution, which does not converge to a zero solution itself.\\

To conclude, we think that the polynomial system of the discrete steady cavity flow problem is challenging for the community of solvers of polynomial systems and numerical algebra. 
Another interesting challenge is to solve the discrete steady cavity flow problem derived by the alternative finite difference discretization of the Jacobian proposed by Arakawa \cite{arakawa}. Our first computational results suggest it is worth further pursuing this alternative.
For its observed and described properties the discrete steady cavity flow problem will be a good test problem to validate new techniques for solving systems of algebraic equations and inequalities. Furthermore, as solving the cavity flow problem for large Reynolds numbers $R$ and velocities $v$ remains an active field of research, we believe that our numerical results may be instructive for audiences in the community of numerical analysis for fluid dynamics to understand fake solutions in partial differential equations.

\section*{Acknowledgements}
We are grateful to Prof. T.~Iwayama for suggesting the alternative finite difference discretization by Arakawa \cite{arakawa} and our discussions about it and to Prof. M.~Noro on a suggestion of using
the RUR method.
We also thank Prof. M.~Kojima for his encouragement and his advice on this study.

\end{document}